\documentclass[11pt]{amsart}

\usepackage[USenglish]{babel}

\usepackage{amsmath,amsthm,amssymb,amscd}
\usepackage{booktabs}
\usepackage{url}

\usepackage{MnSymbol}

\usepackage{tikz}
\usetikzlibrary{arrows,shapes.geometric,positioning,decorations.markings}

\usepackage[mathscr]{eucal}
\usepackage[normalem]{ulem}
\usepackage{latexsym,youngtab}
\usepackage{multirow}
\usepackage{epsfig}
\usepackage{parskip}

\usepackage[margin=2cm]{geometry}
\usepackage{color}

\newtheoremstyle{custom}
  {3pt}
  {3pt}
  {\slshape}
  {}
  {\bfseries}
  {.}
  { }
   {}
\theoremstyle{custom}
\newtheorem{theorem}{Theorem}[section]
\newtheorem{proposition}[theorem]{Proposition}
\newtheorem{observation}[theorem]{Observation}
\newtheorem{proposition/definition}[theorem]{Proposition/Definition}

\newtheorem{conjecture}[theorem]{Conjecture}

\theoremstyle{definition}
\newtheorem{definition}[theorem]{Definition}

\newtheorem{problem}[theorem]{Problem}
\newtheorem{question}[theorem]{Question}

\theoremstyle{remark}
\newtheorem{remark}[theorem]{Remark}

\newtheoremstyle{exercise}
  {3pt}
  {6pt}
  {}
  {}
  {\bfseries}
  {:}
  { }
   {}
\theoremstyle{exercise}
\newtheorem{exercise}[theorem]{Exercise}
\newtheoremstyle{exercises}
  {3pt}
  {6pt}
  {}
  {}
  {\bfseries}
  {:}
  {\newline}
   {}

\theoremstyle{exercise}
\newtheorem{exercises}[theorem]{Exercises}

\newcommand{\exerone}[2][{}]{\begin{exercise}{#1}{#2}\end{exercise}}

\input epsf
\def\boxit#1{\vbox{\hrule height1pt\hbox{\vrule width1pt\kern3pt
  \vbox{\kern3pt#1\kern3pt}\kern3pt\vrule width1pt}\hrule height1pt}}

\def\trank{\text{rank}}

\def\bv{\mathbf{v}}

\def\BC{\mathbb C}\def\BF{\mathbb F}\def\BO{\mathbb O}
\def\BA{\mathbb A}\def\BR{\mathbb R}\def\BH{\mathbb H}
\def\BP{\mathbb P}
\def\pp#1{\mathbb P^{#1}}

\def\fgl{\mathfrak g\mathfrak l}

\def\pp#1{{\mathbb P}^{#1}}
\def\tdim{{\rm dim}}

\def\hd{,...,}
\def\ww{\wedge}

\def\inv{{}^{-1}}

\def\cJ{{\mathcal J}}

\def\cO{{\mathcal O}}

\def\11{\mathbf 1}

\def\fh{{\mathfrak h}}

\def\fsl{{\mathfrak {sl}}}

\def\fg{{\mathfrak g}}

\def\fk{{\mathfrak k}}

\def\a{\alpha}

\def\o{\omega}

\def\b{\beta}
\def\g{\gamma}
\def\s{\sigma}

\def\ot{{\mathord{ \otimes } }}
\def\op{{\mathord{\,\oplus }\,}}

\def\ra{{\mathord{\;\rightarrow\;}}}

\def\dim{{\rm dim}\;}
\def\La#1{\Lambda^{#1}}

\def\frak{\mathfrak}
\def\fgl{\frak g\frak l}\def\fsl{\frak s\frak l}

\def\op{\oplus}
\def\BO{\Bbb O}\def\BA{\Bbb A}\def\BF{\Bbb F}\def\BH{\Bbb H}\def\BZ{\Bbb Z}

\def\bb#1#2#3{b^{#1}_{{#2}{#3}}}

\def\ep{\epsilon}
\def\op{\oplus}

\def\ul{\underline}

\def\s{\sigma}
\def\t{\tau}

\def\a{\alpha}
\def\b{\beta}

\def\g{\gamma}

\def\FS{\mathfrak  S}

\def\ol{\overline}

\def\BP{\mathbb  P}
\def\BC{\mathbb  C}

\def\pp#1{\mathbb  P^{#1}}

\def\tcodim{\text{codim}}
\def\BR{\mathbb  R}

\def\ep{\epsilon}

\def\ci{\mathcal  I}

\def\fg{\mathfrak  g}

\def\hd{, \hdots ,}

\def\inv{{}^{-1}}

\def\La#1{\Lambda^{#1}}

\def\pp#1{\mathbb  P^{#1}}

\def\ur{\underline{\mathbf{R}}}\def\asrk{\uwave{\mathbf{R}}}

\def\ra{\rightarrow}

\def\tdet{\operatorname{det}}

\def\tperm{\operatorname{perm}}
\def\ttrace{\operatorname{trace}}
\def\tend{\operatorname{End}}

\def\tdim{\operatorname{dim}}
\def\tker{\operatorname{ker}}
\def\tlim{\lim}

\def\tmin{\operatorname{min}}
\def\tmax{\operatorname{max}}

\def\trank{\operatorname{rank}}

\def\ww{\wedge}

\def\ctimes{\times \cdots\times}

\def\bbb{{\mathbf{b}}}

\def\be{\begin{equation}}
\def\ene{\end{equation}}
\def\aaa{{\mathbf{a}}}
\def\bbb{{\mathbf{b}}}
\def\ccc{{\mathbf{c}}}

\DeclareMathOperator{\tlog}{log}

\def\trank{\mathbf{R}}

\def\aa#1#2{a^{#1}_{#2}}
\def\bb#1#2{b^{#1}_{#2}}
\def\cc#1#2{c^{#1}_{#2}}

\def\p{\mathbf{P}}
\def\np{\mathbf{N}\mathbf{P}}

\def\G{\Gamma}

\newcommand{\Id}{\operatorname{Id}}

\def\Mn{M_{\langle \nnn \rangle}}\def\Mone{M_{\langle 1\rangle}}\def\Mtwo{M_{\langle 2 \rangle}}\def\Mthree{M_{\langle 3\rangle}}

\def\Mn{M_{\langle \nnn \rangle}}\def\Mone{M_{\langle 1\rangle}}\def\Mthree{M_{\langle 3\rangle}}

\def\Mtwo{M_{\langle 2\rangle}}\def\Mthree{M_{\langle 3\rangle}}

\def\aa#1#2{a^{#1}_{#2}}
\def\bb#1#2{b^{#1}_{#2}}

\def\tinf{{\rm inf}}

\def\trank{{\mathrm {rank}}}

\def\aaa{\mathbf{a}}
\def\bbb{\mathbf{b}}
\def\ccc{\mathbf{c}}

\def\nnn{\mathbf{n}}

\def\bv{\mathbf{v}}
\def\bw{\mathbf{w}}

\def\Det{{\mathcal Det}}\def\Perm{{\mathcal Perm}}

\def\tcap{{\rm cap}}

\def\bt{\bold t}

\def\Mn{M_{\langle \nnn \rangle}}\def\Mone{M_{\langle
1\rangle}}\def\Mthree{M_{\langle 3\rangle}}

\def\Mtwo{M_{\langle 2\rangle}}\def\Mthree{M_{\langle 3\rangle}}

\def\lam{\lambda}

\def\aa#1#2{a^{#1}_{#2}}
\def\bb#1#2{b^{#1}_{#2}}

\def\tinf{{\rm inf}}

\def\trank{{\mathrm {rank}}}

\def\aaa{{\bold a}} \def\ccc{{\bold c}}

\def\nnn{\bold n}

\def\bv{\bold v}\def\bw{\bold w}

\renewcommand{\a}{\alpha}
\renewcommand{\b}{\beta}
\renewcommand{\g}{\gamma}

\renewcommand{\BC}{\mathbb{C}}

 \renewcommand{\tilde}{\widetilde}

\newcommand{\textsum}{{\textstyle\sum}}

\def\bbta{\mathcal B}
\def\uQ{\ul{\bold Q}}\def\uq{\ul{\bold Q}}
\def\assrk{\uwave{\bold Q}}

\def\La#1{\bigwedge^{#1}}

\usepackage[utf8]{inputenc}
\usepackage{amsmath}
\usepackage{amsthm}
\usepackage{amsfonts}
\usepackage{indentfirst}
\usepackage{xcolor}

\usepackage[sort]{cite}

\begin{document}

\author{J. M. Landsberg}

 \email[J.M. Landsberg]{jml@math.tamu.edu}

\title[Geometry and Representation theory in   computer science]{Algebraic Geometry  and Representation theory in
the study of matrix multiplication complexity and other problems in  theoretical computer science}

\thanks{Landsberg supported by NSF grant AF-1814254 }



 \begin{abstract} Many fundamental questions in theoretical
 computer science are naturally expressed as special cases of the following problem: Let $G$
be a complex reductive group,    let $V$ be a $G$-module, and let $v,w$ be elements
of $V$. Determine if $w$ is in the $G$-orbit closure of $v$. I  explain the
computer science problems, the
questions in representation theory and algebraic geometry that they give rise to,
and the new perspectives on old areas such as invariant theory that have arisen in
light of these questions.   I focus primarily on   the complexity of matrix multiplication.
  \end{abstract}

\maketitle

\section{Introduction}

\subsection{Goals of this article}  

To give an overview of uses of algebraic geometry and representation theory in algebraic complexity theory,
with an emphasis on areas that are ripe for further contributions from geometers.

To give a history of, and discussion of recent breakthroughs in, the use of geometry in the study of
the complexity of matrix multiplication, which is   one of the most important problems in algebraic
complexity theory.

To convince the reader of the utility of a geometric perspective  by explaining how the fundamental theorem
of linear algebra is a pathology via secant varieties.

\subsection{The complexity of matrix multiplication}
 
\subsubsection{Strassen's magnificent failure}
In 1968, while attempting to prove the standard row-column way of multiplying matrices is optimal (at least for
$2\times 2$ matrices over $\BF_2$), Strassen instead discovered that $\nnn\times \nnn$ matrices over any field  could be multiplied
using $\cO(\nnn^{2.81})$ arithmetic operations instead of the usual $\cO(\nnn^3)$ in the standard algorithm \cite{Strassen493},
see \S\ref{appen} for his algorithm.
Ever since then it has been a fundamental question to determine just how efficiently matrices can be multiplied.
There is an astounding conjecture that as the size of the matrices grows,  it becomes almost as easy to
multiply matrices as it is to add them.

\subsubsection{The complexity of matrix multiplication as a problem in geometry}
The above-mentioned astounding conjecture may be made precise as follows:
let 
$$\o:=\tinf_{\t}\{ n\times n {\rm \ matrices \ may \ be \ multiplied \ using \ }
\cO(n^\t)\ {\rm arithmetic\  operations} \}
$$
Classically one has $\o\leq 3$ and Strassen showed $\o\leq 2.81$. The astounding
conjecture is that $\o=2$. The quantity $\o$ is called the {\it exponent of matrix multiplication}
and it is a fundamental constant of nature.
As I   explain below, Bini \cite{MR605920}, building on work of Strassen, showed $\o$ can be expressed
in terms of familiar notions in algebraic geometry.

Matrix multiplication is a bilinear map 
$\Mn: \BC^{n^2}\times \BC^{n^2}\ra \BC^{n^2}$, taking a pair of $n\times n$ matrices to
their product:  $(X,Y)\mapsto XY$.
In general, one may view a bilinear map $\b: A^* \times B^*  \ra C $ as
a trilinear form, or {\it tensor}, $T_{\b}\in A\ot B\ot C$. In the case of matrix multiplication, the trilinear form
is $(X,Y,Z)\mapsto \ttrace(XYZ)$.
A basic measure of the complexity of a tensor $T\in A\ot B\ot C$ is its {\it rank}, the smallest
$r$ such that $T=\sum_{j=1}^r e_j\ot f_j\ot g_j$ for some $e_j\in A$, $f_j\in B$, $g_j\in C$. This is   because rank one tensors
correspond to bilinear maps that can be computed by performing one scalar multiplication.
Let $\bold R(T)$ denote the rank of the tensor $T$.
Strassen showed $\bold R(\Mn)=\cO(\nnn^\o)$, so one could 
determine $\o$ by determining the growth of the rank of the matrix multiplication tensor.

\subsubsection{Bini's sleepless nights}
Shortly after Strassen's discovery about multiplying $2\times 2$
matrices with $7$ multiplications instead of $8$, Bini wondered if $2\times 2$ matrices,
where the first matrix had one entry zero, could be multiplied with five multiplications
instead of the usual six. (His motivation was that such an algorithm  would lead to efficient usual matrix multiplication
algorithms for larger matrices.) With his collaborators Lotti and Romani 
\cite{MR592760} they performed a computer
search for a better algorithm, using numerical methods.
There seemed to be a problem with their  code, as each time the program appeared to start to
converge to a rank five decomposition, the coefficients in the terms would blow up.

I had the priviledge of meeting Bini, who told me the story of how he could not sleep at night because
no matter how hard he  checked, he could not find an error in the code and would lie awake  
trying to figure out what was going wrong. Then one evening, he finally realized 
{\it there was no problem with the code}! I explain the geometry of his discovery
in \S\ref{laaside}.

\subsection{Things computer scientists think about}
Computer scientists are not usually  interested in a single  polynomial or fixed size matrix, but rather
{\it sequences} $\{P_n\}$ of polynomials (or matrices) where the number of variables and the degree 
grow with $n$. In what follows I sometimes suppress reference to the sequence.

\subsubsection{Efficient algorithms} This is probably the first thing that comes to mind to a person on the street,
however even here there is a twist:  while  computer scientists are interested in, and very good at, constructing
efficient algorithms, they are also often content just to prove the {\it existence}  of efficient algorithms.
As we will see in \S\ref{lasersect}, this is very much the case for the complexity of matrix multiplication.

\subsubsection{Hay in a haystack}\footnote{Phrase due to H. Karloff} Much of theoretical computer science deals with the problem of
finding explicit examples of objects. A generic (or random) sequence of polynomials will be difficult
to evaluate, but it is a fundamental problem   to find an explicit sequence of polynomials
that is difficult to evaluate.  
(This uses the computer scientist's definition of \lq\lq explicit\rq\rq , which
may not be the first thing that comes to a mathematician's mind.)
Several instances of this problem are discussed in this article:  \S\ref{pnpsect}, \S\ref{nnorsect}, and 
\S\ref{elusect}.

\subsubsection{\lq\lq Lower bounds: Complexity Theory's Waterloo\rq\rq}\footnote{Title of Chapter 14 of \cite{MR2500087}}
Computer scientists are perhaps most interested in proving there are no efficient algorithms for certain
tasks such as the traveling salesperson problem. More precisely, they are interested in distinguishing
tasks admitting efficient algorithms from those that do not. They have not been so successful in this, but
at least they can often show that if task $X$ is difficult, then so is task $Y$, 
e.g., the abundance of $\np$-complete problems.

\subsection{Overview} In \S\ref{laaside} I explain how the fundamental theorem of linear algebra is
a pathology via the geometry of secant varieties, which also explains the geometry underlying Bini's discovery.
In \S\ref{otherproblems} I discuss several problems in complexity theory that are strongly tied to
algebraic geometry and representation theory. Sections  
\S\ref{rrandbr},   \S\ref{retreat}, \S\ref{badnews}, \S\ref{bordersub}, \S\ref{bapolar} and \S\ref{lasersect}
all discuss aspects of the complexity of matrix multiplication, with the first five discussing lower bounds
and the last discussing geometric aspects of upper bounds.
I also include a section, \S\ref{TopenQs}, which  discuses other open questions about tensors, many of which
are also related to complexity theory. For readers not familiar with Strassen's algorithm, I include
it in an appendix \S\ref{appen}. 
The background required for the sections varies considerably,  with almost nothing assumed
in \S\ref{laaside},  where for some of the problems in \S\ref{otherproblems},  I assume basic terminology from
invariant theory.

\subsection{Notation and Conventions} Let $A,B,C,V$ denote complex vector spaces of dimensions
$\aaa,\bbb,\ccc,\bv$. 
Let $a_1 \hd a_{\aaa}$ be a basis of 
  $A$, and $\alpha^1 \hd \alpha^{\aaa}$  its dual basis of $A^*$. 
Similarly $b_1   \hd  b_{\bbb}$ and $c_1   \hd  c_{\ccc}$ are bases of $B$ and $C$ 
respectively, with   dual bases $\beta^1   \hd  \beta^{\bbb}$ and 
$\gamma^1   \hd  \gamma^{\ccc}$. 

The tensor product $A\ot B\ot C$ denotes the space of trilinear maps $A^*\times B^*\times C^*\ra \BC$,
which may also be thought of as the space of linear maps $A^*\ra B\ot C$ etc...
A tensor $T\in A\ot B\ot C$ is {\it $A$-concise}, if the   linear map  $T_A:  A^*\ra B\ot C$
is injective and it is {\it concise} if it is $A$, $B$, and $C$-concise.
  Informally this means that $T$ may not be put in a smaller space.
  
The group of invertible linear maps $A\ra A$ is denoted $GL(A)$ and the set of
all linear maps is denoted $\tend(A)=A^*\ot A$.

Informally, the symmetry group of a tensor $T \in 
A \otimes B \otimes C$ is its 
stabilizer under the natural action
of $GL(A) \times GL(B) \times GL(C)$.
For a tensor $T \in A \otimes B \otimes C$, let $G_T$ denote 
its symmetry group.
 One says $T'$ is {\it isomorphic} to $T$ if  
they are in the same $GL(A) \times GL(B) \times GL(C)$-orbit. I
identify isomorphic tensors.

The transpose of a matrix $M$ is denoted $M^{\bold t}$.

 For a set $S\subset V$, define
the {\it ideal} of $S$, 
$I_S:=\{ {\rm{polys\ }}  P \mid  P(s)=0 \forall s\in S\}$ and define the 
 {\it Zariski closure}  of $S$, 
$\ol{S}:=\{ v\in V\mid P(v)=0 \forall P\in I_S\}$.   
  In our situation the Zariski closure will coincide with the Euclidean closure.  
(This happens whenever $\ol{S}$ is irreducible and $S$ contains a Zariski-open subset of $\ol{S}$,
see, e.g., \cite[Thm 2.33]{MR1344216} or \cite[\S 3.1.6]{MR3729273}).

Projective space is $\BP V=V\backslash\{0\}/\BC^*$.  For the purposes of this article, a {\it projective variety}
is the common zero set of a collection of homogeneous polynomials on $V$ considered as a subset of $\BP V$.

Let $Seg(\BP A\times \BP B\times \BP C)\subset \BP (A\ot B\ot C)$ denote the variety of rank one tensors, 
 called the {\it Segre variety}. 
 
  For a subset $Y\subset \BC^N$,   
let $\langle Y\rangle\subset \BC^N$ denote its linear span and I use the same notation for $Y\subset \BP V$.

For a group $G$,  a $G$-module $V$,  and $v\in V$, $G\cdot v$ denotes the orbit
of $v$,  so $\ol{G\cdot v}$ is its orbit closure.

The space of $d$-way symmetric tensors is denoted $S^dV$, which may be identified with the homogeneous
degree $d$ polynomials on $V^*$. The variety of rank one symmetric tensors
is denoted $v_d(\BP V)\subset \BP S^dV$ and is called the {\it Veronese variety}.

The space of $d$-way skew-symmetric tensors is denoted $\La d V$, and the {\it Grassmannian}
is the variety $G(d,V):=\BP \{ X\in \La d V\mid X=v_1\ww \cdots \ww v_d\mid {\rm some \ } v_1\hd v_d\in V\}$.
It admits the geometric interpretation of the set of $d$-planes through the origin in $V$. 

\subsection{Acknowledgements} I thank D. Alexeev for inviting me to contribute this article,
and J. Grochow, A. Shpilka,  and M.  Forbes for suggestions how to improve the exposition.

\section{The fundamental theorem of linear algebra is an extreme pathology}\label{laaside}
When researchers first encounter tensors they are often surprised  how their intuition from linear
algebra fails and they  view tensors as strange objects. The goal of this section is
to convince the reader that it is not tensors, but matrices that are strange.

\subsection{The fundamental theorem of linear algebra}\label{lasect}

\begin{theorem}[Fundamental theorem of linear algebra] 
 Fix bases $\{a_i\}$, $\{b_j\}$ of $A,B$  
  and for
 $r\leq \tmin\{ \aaa,\bbb \}$,  set $I_r=\sum_{k=1}^r a_k\ot b_k$.
 The following quantities all equal the {\bf rank} of $T\in A\ot B$:  
 \begin{enumerate}
   \item[$(\bold R)$] The smallest $r$ such that $T$ is a sum of $r$  rank one
 elements.    i.e.,   such that  
 $T\in \tend(A)\times \tend(B)\cdot I_r$. 
  \item[$(\ur)$]  The smallest $r$ such that $T$ is a limit of
 a sum of $r$ rank one elements,     i.e.,     such that $T \in \ol{GL(A)\times GL(B) \cdot I_r}$.  
  \item[$(\bold{ml}_A)$]  $\tdim A-\tdim\tker (T_A: A^*\ra B)$. 
  \item[$(\bold{ml}_B)$]  $\tdim B-\tdim \tker( T_B: B^*\ra A) $. 
  \item[$(\ul{\bold Q})$]  The largest  $r$ such that $I_r\in \ol{GL(A)\times GL(B)\cdot T}$. 
  \item[$(\bold Q)$]  The largest $r$ such that $I_r\in \tend(A)\times \tend(B)\cdot T$.  
 \end{enumerate}
 \end{theorem}
 
 \subsection{The fundamental theorem fails miserably for tensors}
 Now  consider a tensor $T\in A\ot B\ot C$.  
   Recall that 
    $T\in A\ot B\ot C$ has {\it rank one} if
there exists  $a\in A$, $b\in B$, and $c\in C$ such that $T=a\ot b\ot c$.
 For
 $r\leq \tmin\{ \aaa,\bbb,\ccc\}$,  write $I_r=\sum_{\ell=1}^r a_\ell\ot b_\ell\ot c_\ell$.

 \begin{definition}\
 
 \begin{enumerate} 
 \item[$( \bold R(T) )$] The    {\it rank} of $T$ is the  smallest $r$ such that $T$ is a sum of $r$ rank one
 tensors    i.e.,   such that
 $T\in \tend(\BC^r)\times \tend(\BC^r)\times \tend(\BC^r)\cdot I_r$, allowing re-embeddings of $T$
 to $\BC^r\ot \BC^r\ot \BC^r$.  
 \item[$( \ur(T) )$]   The    {\it border rank} of $T$ is the smallest $r$ such that $T$ is a limit of rank $r$ tensors, 
 i.e.  such that $T \in \ol{GL(\BC^r)\times GL(\BC^r)\times GL(\BC^r)\cdot I_r}$, allowing re-embeddings.  
  \item[$( \bold{ml} )$]  The {\it multi-linear ranks}
   of $T$ are  $(\bold{ml}_A(T),\bold{ml}_B(T),\bold{ml}_C(T)):=  (\trank T_A, \trank T_B, \trank T_C)$.
\item[$( \ul{\bold Q}(T)) $]  The    {\it border subrank} of $T$ is the    largest  $r$ such that $I_r\in \ol{GL(A)\times GL(B)\times GL(C)\cdot T}$.    
 \item[$( \bold Q(T) )$] The  {\it subrank} of $T$ is the    largest $r$ such that $I_r\in \tend(A)\times \tend(B)\times \tend(C)\cdot T$.
 \end{enumerate}
\end{definition}

\begin{proposition} One has  
 $$\bold Q(T)\leq \uQ(T)\leq \tmin\{\bold{ml}_A(T),\bold{ml}_B(T),\bold{ml}_C(T)\}
 $$
 $$
 \leq \tmax\{\bold{ml}_A(T),\bold{ml}_B(T),\bold{ml}_C(T)\}\leq \ur(T)\leq \bold R(T)
 $$
and    all inequalities  may be strict, even when $\aaa=\bbb=\ccc$. 
\end{proposition}

     \smallskip
 
  Say $\aaa=\bbb=\ccc=m$, then it has been known for some time that if 
  $T$ is generic then  $\ur(T)=\bold R(T)\simeq \frac {m^2}3$,  
 and this is largest possible $\ur$.  However the precise generic values were not
 determined until Lickteig determined them in 1985 \cite{MR87f:15017}.
 The  symmetric case  was studied by
 Terracini  in 1916,    who mostly solved it but it was not finished (for polynomials
 of arbitrary degree) until 1995 when it was solved by
 Alexander and Hirschowitz \cite{AH}. 
 
 \newpage
 
 \subsection{Geometry of Bini's insight}
 Consider the following three pictures:

\begin{figure}[!htb]\begin{center}
\includegraphics[scale=.3]{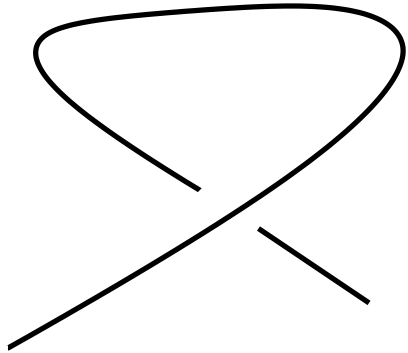}
\caption{Imagine this curve represents the $3m-2$ dimensional set of tensors of rank one sitting in the $m^3$ dimensional
space of tensors. }
\end{center}
\end{figure}

\begin{figure}[!htb]\begin{center}
\includegraphics[scale=.3]{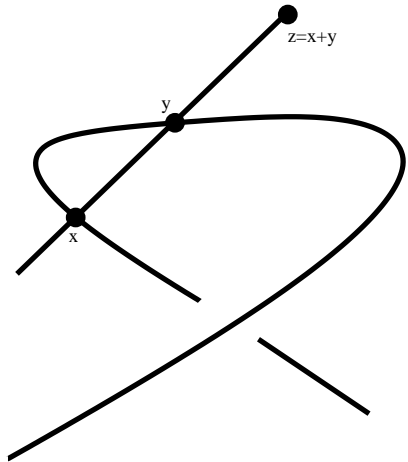}
\caption{Tensors of rank two correspond to points on a secant line to the set of tensors of rank one}
\end{center}
\end{figure}

\begin{figure}[!htb]\begin{center}
\includegraphics[scale=.3]{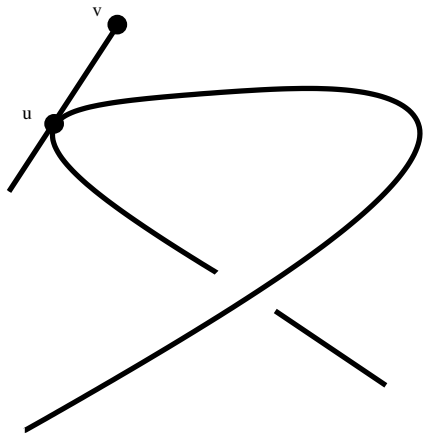}
\caption{The limit of secant lines is a tangent line!  }
\end{center}
\end{figure}

Given a  curve or other variety with large codimension,
most points that lie on a secant line  lie   on just one secant line and  
  points that lie on a secant line do  not in general lie on a tangent line.
Contrast this with  the case of a plane curve, where all points lie on a family
of secant lines.

Bini's (re)discovery was that tensor rank is not semi-continuous precisely because
secant lines may limit to tangent lines, and he coined the term
border rank to include the limits. (The classical Italian algebraic geometers
knew of the lack of semi-contiuity 100 years ago.) Bini then went on to prove that border rank is 
also a legitimate measure of the complexity of matrix multiplication, namely
$\ur(\Mn)=\cO(n^\o)$ \cite{MR605920}. 
 
 Now we see that  Strassen's result   $\bold R(\Mn)=\cO(n^\o)$ is  not immediately  useful for algebraic geometry:
if  $P$ is a polynomial,  $P(T_t)=0$ for $t>0$  implies $P(T_0)=0$, but
the limit of tensors of rank at most $r$ need not have rank at most $r$.  
I.e., one cannot describe  rank via zero sets of  polynomials.  
In contrast, for matrices matrix   rank equals matrix border rank and
is  given by polynomials.  For this reason, we will be primarily concerned with border rank
of tensors, which by definition is closed under taking limits.

As I explain below,   essentially  all known examples of  smooth geometric objects with large codimension
with the property that any time a point is on a secant line to the object, it is on a family of secant lines,
comes from \lq\lq minimal rank\rq\rq\ matrices in some (restricted) space of matrices.
 
\subsection{Geometric context: secant varieties}
Let $X\subset \BP V$ be a projective variety. Define 
$$
\s_r(X):=\ol{ \{z\in \BP V \mid \exists x_1\hd x_r\in X \mid z\in \langle x_1\hd x_r\rangle  \} }
$$
the variety of secant $\pp{r-1}$'s to $X$.
A na\"\i ve dimension count shows that we expect $\tdim(\s_r(X))$ to be $\tmin\{r\tdim X+r-1, \tdim \BP V\}$,
because we get to pick $r$ points on $X$ and a point in their span. If this count fails one says the variety
is {\it defective}. Defectivity is a pathology.

The fact that the fundamental lemma of linear algebra is a pathology may be rephrased as saying the
Segre variety $Seg(\BP A\times \BP B)\subset \BP (A\ot B)$ of rank one matrices has
defective secant varieties. The secant varieties of three and higher order tensors are generally
not defective, the $(m,m,m)$ case mentioned above proved  by Lickteig has only one exception, $m=3$ and $r=4$.
See \cite{MR2452824} for the state of the art in the general tensor case.

To my knowledge, essentially all known smooth varieties with degenerate $r=2$ secant varieties
come from matrices: the rank one matrices $Seg(\BP A\times \BP B)\subset \BP (A\ot B)$,
the rank one symmetric matrices $v_2(A)\subset \BP (S^2A)$, the rank two skew-symmetric
matrices (i.e., minimal rank skew-symmetric matrices)  $G(2,A)\subset \BP( \La 2 A)$, the Cayley plane discussed below, and the
closed orbit in the adjoint representation of a complex simple Lie algebra $G/P_{\tilde \a}\subset \BP \fg$
\cite{MR1663641} (which consists of $v_2(\BP A)$,  the  two planes isotropic for 
a quadratic form on $A$, denoted  $G_Q(2,A)$, the traceless rank one matrices when $\fg=\fsl_n$,
 and five
exceptional cases of \lq\lq minimal rank\rq\rq\ matrices in some space of matrices with extra structure).
When the dimension of the secant variety differs from the expected dimension by more than one,
one can take hyperplane sections to get further examples.

Secant varieties of projective varieties have been studied for a long time, dating back
at least to the Italian school in the 1800's. They show up in numerous geometric situations.
For example, Zak (see \cite{Zak}) proved a linear approximation to Hartshorne's famous conjecture on 
complete intersections \cite{MR0384816} that was also conjectured by Hartshorne, which 
amounts to proving that if $X$ is smooth of dimension $n$, not contained in a hyperplane,   and the codimension of $X$ is
less than  $\frac n2$, then $\s_2(X)=\BP V$. Moreover, Zak classified the exceptional cases
with codimension $\frac n2$, now called
Severi varieties. They turn out to be the projective planes over the composition algebras $\BA\pp 2\subset
\BP \cJ_3(\BA)$, where $\BA$ is one of $\BR^{\BC}=\BC$, $\BC^\BC=\BC\op \BC$, $\BH^\BC$, or $\BO^{\BC}$,
the last two cases are the complexified quaternions and octonions, and $\cJ_3(\BA)$
denotes the (complex) vector space of $3\times 3$ $\BA$-Hermitian matrices. 
It is an open problem to determine if a gap larger than $8$ (which occurs for $\BO^{\BC}\pp 2$)  between the expected
and actual dimension occurs for a smooth variety, see \cite{LVdV}.  Moreover, using 
the geometry of secant and tangential
varieties, one can obtain a proof of the Killing-Cartan classification of complex simple Lie algebras
via a constructive procedure starting with $\pp 1$, see \cite{MR1890196}.

\section{Problems in complexity with interesting  geometry}\label{otherproblems}
The bulk of this article discusses the history and recent developments using geometry in the study
of matrix multiplication. Here I discuss several other questions with interesting geometry. I ignore
numerous  fundamental developments that have not yet been
cast in geometric language.  I briefly mention two such developments now:  Hardness v. randomness (see, e.g., \cite{MR1957568}),
  shows that if certain   conjecturally hard problems are truly hard, then many other problems, for
which only a randomized efficient algorithm is known,  admit deterministic efficient algorithms,
and conversely, if these problems cannot be derandomized, then the conjectured hard problem
is not hard after all.
The PCP theorem \cite{MR1639346}  says that any putative  proof may be rewritten in such a way
that its correctness is checkable by looking at only a few probabilistically chosen symbols,
 See \cite{MR2500087} for excellent discussions of these and other  omissions.

\subsection{P v. NP and variants}\label{pnpsect}
The famous $\p\neq\np$ conjecture of Cook, Karp,  and Levin
has origins in the 1950's  work of John Nash (see \cite[Chap. 1]{nash2016open}), researchers in the Soviet Union (see, e.g,
 \cite{MR763733}), and, most poetically, 
G{\" o}del, who asked if one could quantify the idea of intuition (see \cite[Appendix]{Sipser}).
One of L. Valiant's algebraic versions of the problem \cite{vali:79-3} is as follows: Valiant showed  that any polynomial $p(x_1\hd x_N)$
may be realized as the determinant of some $n\times n$ matrix of affine linear forms in the $x_i$, where the size $n$ depends on $p$.
He showed the size is a good measure of the complexity of a polynomial.
The \lq\lq permanent v. determinant\rq\rq \ version of the $\p\neq\np$ conjecture  is that the size $n(m)$ of the matrix  needed to compute the
permanent $\tperm_m\in S^m\BC^{m^2}$ of an $m\times m$ matrix\footnote{The permanent is the  polynomial  that is the same as the determinant, only without the minus signs.} grows faster than any polynomial in $m$. For those familiar with the traveling
salesman version of $\p\neq\np$, the \lq\lq easy to verify\rq\rq\ proposed answer is replaced by an \lq\lq easy to write down\rq\rq\ 
polynomial sequence, and the conjecture amounts to saying a polynomial sequence  that is easy to write down in general will not be easy to
compute.  Here I   discuss the permanent v. determinant version of the problem.

In algebraic geometry one generally prefers to work with homogeneous polynomials, so instead
of asking to write the permanent as a determinant of affine linear functions, following K. 
Mulmuley and M. Sohoni, one adds a variable
to homogenize the problem: 
\begin{conjecture}\label{valrephrase} [Rephrasing of Valiant's conjecture as in \cite{MS1}]
Let $\ell$ be a linear coordinate on $\BC^1$ and consider any linear inclusion
$\BC^1\oplus \BC^{m^2}\ra  \BC^{n^2}$, so in particular $\ell^{n-m}\tperm_m\in S^n\BC^{n^2}$.
Let $n(m)$ be a polynomial. Then for all sufficiently large $m$, 
$$
[\ell^{n-m}\tperm_m]\not\in  \tend(\BC^{n^2}) \cdot [\tdet_{n(m)}] .
$$
\end{conjecture}

\begin{remark} The best known $n(m)$ is $n=2^{m}-1$ due to B. Grenet \cite{Gre11}. His expression has
interesting geometry. Note that given an   expression for the permanent as a determinant, one gets a family of such
by the action of $G_{\tperm_m}\times G_{\tdet_n}$ on $(\BC^1\oplus \BC^{m^2})^*\ot    \BC^{n^2}$,
and we can consider the stabilizer $\G\subset G_{\tperm_m}\times G_{\tdet_n}$ of a decomposition.
One has $G_{\tperm_m}=[(\BC^*)^{\times 2m-2}\times \FS_m\times \FS_m]\rtimes \BZ_2$ and this is the largest possible
stabilizer. One has $\G_{Grenet}=(\BC^*)^{\times  m}\times \FS_m$ \cite{MR3724217}
and moreover, if one insists on a stabilizer this size, N. Ressayre and I  also  showed Grenet's $n(m)=2^m-1$ is the  smallest possible \cite{MR3724217}. In particular, one could prove Valiant's conjecture
by proving that for any determinantal expression for the permanent, that there exists a slightly
larger one with symmetry $\G_{Grenet}$.
Similarly, one can also have $\G=G_{\tperm_m}$ for a similar (but larger) cost in size.
\end{remark}

To further geometrize the problem, Mulmuley and Sohoni proposed a stronger conjecture:
Let 
$$\Det_n:=\ol{GL_{n^2}\cdot [\tdet_n]},
$$ 
and let 
$$\Perm^m_n:=\ol{GL_{n^2}\cdot [\ell^{n-m}\tperm_m]}.
$$
\begin{conjecture}\label{msmainconjb}  \cite{MS1}  Let $n(m)$ be a polynomial. Then for all sufficiently large $m$, 
$$
\Perm^m_{n(m)}\not\subset \Det_{n(m)}.
$$
\end{conjecture}

This is stronger as in general, for a polynomial $P\in S^d\BC^N$, the inclusion  $\tend(\BC^N)\cdot P\subset  \ol{GL_N\cdot P}$
is   strict. 

Both $\Perm^m_n$ and $\Det_n$  are {\it invariant} under $GL_{n^2}$ so their
ideals are  $GL_{n^2}$-modules. The original idea of \cite{MS1,MS2}
was to solve the problem by finding a  sequence, depending on $n$, of   $GL_{n^2}$-modules
$M_n$ such that 
$M_n\subset I[\Det_n]$ and $M_n\not\subset I[\Perm^m_n]$.
The initial idea in \cite{MS1} was to look for not any module, but an invariant to separate the orbit closures.
The mathematical issues raised by this program were analyzed in \cite{MR2861717}.
Work of  Ikenmeyer and Panova \cite{MR3695867} and 
  B{\"{u}}rgisser, Ikenmeyer and
Panova \cite{MR3868002} shows this is not possible but other paths using representation
theory are still open, see \cite{MS3,MS4,MS5,MS6,MS7,MS8}.

This program has inspired a tremendous amount of work and breathed new life into invariant theory.
I will loosely refer to the inspired work as {\it geometric complexity theory} (GCT) even for problems 
not directly related to Valiant's conjecture.

\subsection{Explicit Noether Normalization}\label{nnorsect}
The classical Noether normalization lemma may be stated geometrically as follows: given an affine
variety $X^n\subset \BA^{n+a}=V$, there exists a $W=\BA^{a}\subset \BA^{n+a}$ such that the projection of $X$
 to $V/W =\BA^n$  will be a finite to one surjection. In fact a general or \lq\lq random\rq\rq\  $\BA^{a}$ will do.
The problem of {\it explicit Noether normalization} is to find an explicit $\BA^{a}$ that
has this property. In \cite{MS5}, Mulmuley proposes this problem as a possible path to
resolve Valiant's conjecture.

Algebraically, given a ring $R$, one is interested in finding an explicitly generated
(in the sense of computer science, in particular, efficiently describable, and I remind the reader that
one is really dealing with sequences of rings)  subring
$S$ such that $R$ is integral over $S$. Of particular interest is the case $R$ is
the ring of invariants in the coordinate ring of some $G$-variety $Z$ (i.e., the $X$ in question
  is the GIT quotient of $Z$). 
M. Forbes and A. Shpilka \cite{MR3126552} give an effective solution to this in the case
$R$ is the ring of invariants of the set
of $r$-tuples of $n\times n$ matrices under the action of $GL_n$ by simultaneous
conjugation. This ring is generated by traces of words in the matrices
\cite{MR0506414,MR419491}.
A key ingredient here is the use of {\it pseudo-randomness}: deterministic processes that 
are not too expensive yet are able to imitate randomness sufficiently to resolve the problem.

There is work still to be done here.  One wants to find a  smallest  subring $S$ that works.  One way to measure 
smallness  is   the number of separating invariants that are used to generate it.  Ideally one would get a set of generators whose size matches the dimension of the variety in question, and generically this is possible. Here the dimension   is polynomial in  $(n,r)$, one writes \lq\lq $poly(n,r)$\rq\rq.  The articles \cite{MR0506414,MR419491}
give an exponential in $n,r$  size $S$ and  \cite{MR3126552} gives an explicit set of size $poly(n,r)^{\cO(\log(nr))}$
and the open problem is to  reduce it  to polynomial size.

\subsection{Algorithms in invariant theory}\label{algsinvarsect}
Thanks to the GCT program, computer scientists have recognized that the orbit closure containment
problem, and the related orbit closure intersection and orbit equality problems provide a natural
framework for many questions in computer science.
One  insight of \cite{Burgisser_2019},  by P. Burgisser, C.  Franks, A.  Garg, R.  Oliveira, M.  Walter,  and A. Wigderson,    is that many problems in complexity and
other applications may be phrased
as follows: Let $G$ be a complex  reductive algebraic group, let $K\subset G$ be a maximal compact subgroup with Lie algebra $\fk$,
so $\fg$ is a complex Lie algebra and $\fk$ is a real Lie algebra. Then $\fg=\fk\op i\fk$, and $G$-modules
$V$ come equipped with a $K$-invariant inner product, in particular a norm.
Then given $v\in V$, the problem is
to compute its {\it capacity}
$$
\tcap (v):=\tmin_{w\in \ol{G\cdot v}} ||w||.
$$ 
Readers familiar with the Kempf-Ness theorem \cite{MR555701} will recognize $w$ as a vector
satisfying $\mu(w)=0$, where $\mu: \BP V\ra (\fg/\fk)^*$ is the {\it moment map}, the  normalized derivative of the norm
map $G\ra \BR$, $g\mapsto ||g\cdot v||^2$ at the identity.

They point out that already when $G$ is abelian, this encompasses all linear programming problems.
Incidences of the general problem  have been known for some time in the community of geometers, e.g., the famous
Horn problem, see, e.g.,  \cite{MR1685640}. A very special case is null cone membership: determine
if $v$ has capacity zero.

The problem with just checking the boundary of the moment polytope is that in general, the moment polytope
is defined by too many inequalities to make such a check efficient. 

Their new ingredient to this well-studied problem is the introduction of algorithms to
solve, or approximately solve the problem. They use gradient descent methods. In essence,
given $\ep>0$, the algorithms, after some specified amount of computation, either take one within
$\ep$ of a $w$ achieving the minimum, or stay a specified distance from
it.  A  variant that they also study  is an algorithm to   compute an element of $G$ that takes one close to $w$.

 \subsection{Elusive functions}\label{elusect}
R. Raz defines the following  \lq\lq hay in a haystack\rq\rq\ approach to Valiant's conjecture.
Consider a  linear projection of a Veronese variety $v_r(\pp{s-1})\subset \BP S^r\BC^s$
via  $proj: \BP S^r\BC^s \dashrightarrow \pp m$,
and let $\G_{r,s}:= proj\circ v_r: \pp {s-1}\dashrightarrow \pp m$ be the composition of the projection with the Veronese map.
A map $f: \pp{n}\ra \pp m$ is said to be $(r,s)$-{\it elusive} if $f(\pp n)$ is not contained in the
image of any such $\G_{r,s}$.

\begin{theorem} \cite{MR2719753} Let $m$ be super-polynomial in $n$, and $s\geq m^{\frac 9{10}}$. If there exists an
explicit  $(s,2)$-elusive $f: \pp{n}\ra \pp m$, then Valiant's conjecture is true.
\end{theorem}

\begin{theorem} \cite{MR2719753} Let $r(n)=\tlog(\tlog(n))$, $s(n)=n^{\tlog(\tlog(\tlog(n)))}$,   $m=n^r$,
and let $C$ be a constant.
If there exists an explicit $(s,r)$-elusive $f: \pp{n}\ra \pp m$, then  Valiant's conjecture is true.
\end{theorem}

 By a dimension count, a general  polynomial  in either range will be elusive.
The problem of finding explicit elusive functions seems to be worth further study with the tools of algebraic geometry.
See \cite{guo2021variety} for recent work in this direction. 

\subsection{Shallow circuits}\label{shallowsect}
The usual model of computation in algebraic complexity theory is the {\it arithmetic circuit}, which is a directed
graph used to encode a polynomial allowing additions and multiplications and the size of the circuit
(essentially the number of additions and multiplications) measures its complexity.  There have been
results regarding restricted circuits, where, e.g., one is just allowed a round of sums, then a round of
products then a third round of sums, called  \lq\lq $\Sigma\Pi\Sigma$ circuits\rq\rq . Then if one can prove a strong lower
bound for computing the permanent $\tperm_m$  in such a restricted model (roughly exponential in $\sqrt{m}\tlog^{\frac 32}m$),
one can prove Valiant's conjecture. See, e.g., \cite{MR0660280,DBLP:journals/eccc/GuptaKKS13,MR3303254,koirand4,AgrawalVinay}.

The essential point here for geometry is that  the Zariski closure of   the set of polynomials  that $\Sigma\Pi\Sigma$ circuits
of fixed size and \lq\lq fainin\rq\rq \ 
can produce corresponds to a well-studied, albeit  little understood,  object in representation theory and
algebraic geometry: secant varieties of the (simplest)  {\it Chow variety}.

The Chow variety $Ch_n(W)\subset \BP S^nW$ is the set of polynomials that
are products of linear forms. When $\bw=\tdim W\geq n$, it is the orbit closure
$\ol{GL(W)\cdot [x_1\cdots x_n]}$.  It is a longstanding open problem to understand
its ideal dating back to  Brill,  Gordan, Hermite and especially Hadamard \cite{MR1504330}.  See \cite{MR983608,MR1243152,MR1601139,MR3314828} for its connections
to interesting questions in algebraic geometry and representation theory. See \cite[Chap. 9]{MR3729273} for an exposition
of the state of the art.
  
 Other models of shallow circuits give rise to other interesting secant varieties. 
 
\subsection{Identifiability}\label{identsect} While the following problem in applications initially comes from engineering,
it is similar enough in spirit that I include it here.
A basic geometric fact is that under a Veronese re-embedding points become more independent.
A collection of $r$ points is in {\it general linear position} if for all $k$,  no subset of $k$ of them lies on a $\pp{k-2}$.
More generally,  a collection of $r$ points is in {\it $d$-general position} if for all $k$, no subset of $k$ them lies on a higher dimensional
space of hypersurfaces of degree $d$ than holds for a generic set of $k$ points. Under the $d$-th Veronese
re-embedding $v_{d }: \BP V\ra \BP S^{d }V$, given by $[x]\mapsto [x^d]$, collections of points become in $(d-1)$-general position.
P. Comon \cite{Como94:SP} exploited this in signal processing. See \cite[Chap. 12]{MR2865915}
for a geometric discussion of how this was exploited.

The general question that arises in applications is as follows:  given a tensor    of rank $r$ with a rank $r$ decomposition, determine if
the decomposition is unique (up to trivialities). The first, and still most important, result in this regard is
Kruskal's theorem \cite{MR0444690}, which assures uniqueness in a certain range if the points giving rise to the decomposition
are in general linear position. In \cite{MR2996363} it was shown that Kruskal's theorem is sharp.
On the other hand, tensors are known to be generally identifiable well beyond the Kruskal range \cite{MR3023462,MR3624400}. In a series of papers,  generalizations of Kruskal's theorem that exploit more subtle geometric information
have extended identifiability. For example, a first step beyond Kruskal's bound is obtained in \cite{MR3819100} by exploiting
Castenuovo's theorem that if a set of $2n+3$ points lies on  a $\binom n2$-dimensional space of quadrics, then
the points all lie on a rational normal curve. The current state of the art is \cite{angelini2021description},
where further advanced tools (minimal free resolution, liason, Hilbert functions) are used.

\section{Representation theory and border rank}\label{rrandbr}
When I was first introduced to the problem of matrix multiplication,
$\ur(\Mn)$ was not known in any case  except the trivial $\nnn=1$,  and it
was an often stated open problem just to determine $\ur(\Mtwo)$.
I now explain why I thought this would be easy to resolve using representation theory.

\begin{definition}
Given a variety $X\subset \BP V$, one says  $X$ is a {\it $G$-variety}
if it is   invariant  under the action of some 
group $G\subset GL(V)$, i.e.,  $\forall g\in G$, $\forall x\in X$, $g\cdot x\in X$.
\end{definition}

For example, 
$\s_r(Seg(\BP A\times \BP B\times \BP C)):=\{ [T]\mid \ur(T)\leq r\}$, is
a $G=GL(A)\times GL(B)\times GL(C)\subset GL(A\ot B\ot C)$  variety.

As mentioned above, a fundamental observation about $G$-varieties is that their ideals are
$G$-modules.

In particular, if $G$ is reductive, one  could in principle determine  the ideal 
$I(X)$ in any given degree $d$ by decomposing $S^dV^*$ as a $G$-module
and then testing highest weight vectors on random points of $X$. If the polynomial
vanishes at a general point, then the module is in the ideal and otherwise it is not.

In the special case where $X=G/P\subset \BP V_{\lam}$ is homogeneous (as with the Segre variety),
Kostant showed that  the
ideal is generated in degree two by  $V_{2\lam}^\perp \subset S^2V^*$, where $\lam$ is the highest
weight of the irreducible $G$-module $V$ (Kostant's proof 
appeared in an appendix to  the unpublished \cite{garfinkle}. See \cite[\S 16.2]{MR2865915}
for a proof  or \cite{MR1923198} for a proof with an 
extension  to the infinite dimensional case).

For example, $S^2(A^*\ot B^*)= S^2A^*\ot S^2B^*\op \La 2 A^*\ot \La 2 B^*$,
and $I(Seg(\BP A\times \BP B))$ is generated in degree two by $\La 2 A^*\ot \La 2 B^*$
which spans the two by two minors.
More generally $I(\s_r(Seg(\BP A\times \BP B)))$ is generated in degree $r+1$
by the size $r+1$ minors $\La{r+1}A^*\ot \La{r+1}B^*$.

Tensors of border rank at most two are the zero set of degree three polynomials  \cite{LMsec}, so
I was optimistic.
After all, to decide the border rank of $\Mtwo$, one just needs to determine polynomials
in the ideal of $\s_6(Seg(\pp 3\times \pp 3\times \pp 3))$ and to test them on $\Mtwo$. (Strassen
had previously proved the border rank was at least six and $\s_7(Seg(\pp 3\times \pp 3\times \pp 3))=
\BP (\BC^4\ot \BC^4\ot \BC^4)$.)
 With L. Manivel, we carried out a systematic search. Unfortunately we found:

\begin{theorem}\cite{LMsec} The ideal of $\s_6(Seg(\pp 3\times \pp 3\times \pp 3))$ is empty in degrees less than $12$.
\end{theorem}

It was probably this proposition that ended Manivel's activity in the area. Later,  J. Hauenstein and C. Ikenmeyer  extended the
result:
  
  \begin{theorem}\cite{MR3171099} The ideal of $\s_6(Seg(\pp 3\times \pp 3\times \pp 3))$ is empty in degrees less than $19$.
\end{theorem}

Fortunately, at the same time we showed:
\begin{theorem}\cite{MR3171099} A copy of the degree $19$ module
$S_{5554}A^*\ot S_{5554}B^*\ot S_{5554}C^*$  is in the ideal of $\s_6(Seg(\pp 3\times \pp 3\times \pp 3))$.
\end{theorem}

Here I use partition notation to describe modules for the general linear group. Unfortunately, the polynomials
were too complicated to test symbolically on matrix multiplication, so it was only useful to get a numerical proof.
However
  there is an invariant in degree $20$ that was easier to work with which provided the first algebraic proof
of

\begin{theorem} \cite{MR2188132,MR3171099} $\ur(\Mtwo)=7$.
\end{theorem}

The original proof in  \cite{MR2188132} was obtained using differential-geometric methods. 

Clearly to deal with larger $\Mn$, different methods were needed.

\section{Retreat to linear algebra}\label{retreat}
 A concise tensor   $T\in \BC^m\ot \BC^m\ot \BC^m$ must have border rank at least $m$, and if equality holds,
 one says  $T$ has {\it minimal border rank}. As explained in \S\ref{lasersect}  below, minimal border rank tensors are important for 
Strassen's laser method.

Strassen's proof that $\ur(\Mtwo)\geq 6$ was obtained by taking advantage of the following
correspondence: given a concise $T\in A\ot B\ot C$, one obtains
an $\aaa$-dimensional  space of matrices $T(A^*)\subset B\ot C$, and $T$ may be recovered up to
isomorphism from this subspace. I.e., we have a correspondence

{\it $A$-concise tensors $T\in A\ot B\ot C$ up to $GL(A)\times GL(B)\times GL(C)$-isomorphism $\leftrightarrow$
$\aaa$-dimensional subspaces $U\subset B\ot C$ up to $GL(B)\times GL(C)$-isomorphism.}

While spaces of linear maps are nice, even better are spaces of endomorphsims: let $ \bbb=\ccc=m$
and assume there exists $\a\in A^*$ with $T(\a): B^*\ra C$ of  full rank $m$.
Then  $T(A^*)T(\a)\inv\subset \tend(C)$, and one can recover $T$ up to isomorphism from this space as well.

Now if $\aaa=m$ and  $\bold R(T)=m$,  then $T(A^*)T(\a)\inv$ is an $m$-dimensional  space of simultaneously diagonalizable endomorphisms,
which, in  a good basis   is $\langle \g^1\ot c_1,\hd \g^m\ot c_m\rangle$.
So if $T$ is of minimal border rank, then $T(A^*)T(\a)\inv$ is a
limit (in the Grassmannian $G(m, \tend(C))$) of    spaces of simultaneously diagonalizable endomorphisms.
So the problem to determine if $T$ has minimal border rank is reduced to determining
if $T(A^*)T(\a)\inv$ is such a limit.
Good news: this problem was studied classically in the linear algebra literature (e.g., \cite{MR0132079}).
Bad news: it is still open!

Nonetheless, it is easy to obtain necessary conditions: simultaneously diagonalizable matrices commute,
and commutivaty  is a Zariski closed condition. Call the vanishing of the commutators
{\it Strassen's equations for minimal border rank}.  Moreover, Strassen showed that the failure of commutativity
(i.e., the rank of the commutator) lower bounds the border rank, which enabled
the first lower bound on $\ur(\Mn)$:

\begin{theorem} \cite{Strassen505} $\ur(\Mn)\geq \frac 32 n^2$.
\end{theorem}
Lickteig \cite{MR86c:68040} was able to improve the error term to obtain $\ur(\Mn)\geq \frac 32 n^2+\frac n2-1$.
Then from 1985 to 2012 there was no further progress on the general case.

Taking a more abstract view of Strassen's theorem, he found equations by embedding
$A\ot B\ot C$ into a space of matrices and then took minors. This is a longstanding trick
in algebraic geometry: obtain   determinantal equations of varieties by looking at rank loci
of maps between vector bundles on projective space. In the situation of $G$-varieties,
one can give an elementary description of the most useful embeddings:

\begin{observation} \cite{MR3081636} Given a $G$-variety $X\subset \BP V_{\lam}$, say
$V_{\lam}$ occurs as a submodule of $V_\mu\ot V_\nu$. Then if a general point of
$X$ maps to a rank $t$ element (as a matrix) of $V_\mu\ot V_\nu$, then
the size $rt+1$ minors restricted to $V_{\lam}$ are in the ideal of  $\s_r(X)$.
\end{observation}

Call such an inclusion $V_\lam \ra V_\mu\ot V_\nu$ a {\it Young flattening}. Several
such were useful for the case $X$ is a Veronese variety \cite{MR3081636}. For three way tensors,
the most useful inclusions have been, for values of $p<\aaa/2$,  $A\ot B\ot C\ra (\La p A^*\ot B)\ot (\La{p+1}A\ot C)$,
which we call {\it Koszul flattenings}. To implement them one generally restricts $T$ to 
a $2p+1$ dimensional subspace of $A$, and it is often  an art to find an explicit useful
such subspace.
Using a judiciously chosen $SL_2\subset SL_A$ to define a good restriction, G. Ottaviani and I were able to show: 

\begin{theorem}\cite{MR3376667}  $\ur(\Mn)\geq 2n^2-n$.
\end{theorem}

\section{Detour: Additional  open questions regarding tensors}\label{TopenQs}

 Unlike the case of linear maps, we are rather ignorant of tensor rank:
 
\begin{question} For $T\in \BC^m\ot \BC^m\ot \BC^m$, what is the  largest possible $\bold R(T)$? 
\end{question}
For the state of art, see \cite{MR3769376} and \cite{MR3368091}.

\begin{problem}\cite[Problem 15.2]{BCS}\label{minbr}  Classify concise tensors of minimal border rank.
\end{problem}

We are rather ignorant here as well: the state of
the art is  $m=4$ (Friedland \cite{MR2996364}).

As innocent as Problem \ref{minbr} may sound,   even a special case of it amounts to
characterizing when  a zero dimensional scheme is smoothable, a notoriously difficult problem
that is known only for very small values of $m$: when $m\leq 7$ all are, and  when $m=8$ the problem
is solved in  \cite{MR2579394}. For larger $m$ little is known, see \cite{MR3999690} for the state
of the art as of this writing. I now explain this connection.

Let $A=B=C=\BC^m$. 
Call a tensor $1_A$-generic if $T(A^*)\subset B\ot C$ contains an element of full rank $m$.
Call a tensor $1_*$-generic if it is at least one of $1_A$, $1_B$ or $1_C$-generic, and
{\it binding} if the property  holds in at least two directions, and $1$-generic if it holds in all three.

M. Bl\"aser and V. Lysikov \cite{MR3578455} showed that if a tensor is binding, then it is the structure
tensor of some (not necessarily commutative) algebra with unit. Strassen's commutivity
equations for minimal border rank described above   are a necessary condition for minimal border rank and satisfying
them implies  that the algebra is commutative, which implies  the algebra is of the form $\BC[x_1\hd x_n]/\ci$, where
$\ci$ is some ideal and $\BC[x_1\hd x_n]/\ci$ is an  $m$-dimensional vector space over $\BC$.
This leads to the Hilbert scheme for zero dimensional schemes of length $m$, which
parametrizes such objects, see \S\ref{hs} for more details. Then the question becomes whether the given algebra lies
in the same component as the algebra $(\BC[x]/(x^2))^{\op m}$ corresponding to a 
concise rank $m$
tensor.

If $\bbb=\ccc=m$ and the tensor is only $1_A$-generic, there is still a geometric object one can utilize, the
{\it Quot scheme} parametrizing modules: the 
Atiyah-Hitchin-Drinfeld-Manin (ADHM) correspondence
associates  to the tuple $(x_1\hd x_\aaa)$ of commuting $m\times m$ matrices one gets from 
Strassen's commutivity equations,  
a  $\BC[y_1\hd y_\aaa]$-module structure on $V=\BC^m$ defined by $y_i(v):=x_i(v)$.  
The
  {\it Quot scheme of points} is a moduli space for such modules.  
(The special case where the  module has one generator is the Hilbert scheme of points.) Jelisiejew and   Sivic
\cite{jelisiejew2021components}  use this correspondence to prove new results about each. In 
particular, they classify all components when  $m\leq 7$. 
In recent work with A. Pal and  J. Jelisiejew, we  use  this to solve problem \ref{minbr}
under the assumption of $1_*$-genericity up to $m\leq 6$.

If a tensor fails to be $1_*$-generic, one is led to the   problem of
characterizing spaces of matrices of bounded rank - a classical but difficult topic that
has only be solved for ranks up to $3$ \cite{MR695915}. Here one is in a slightly better situation,
as there are three such spaces to be considered.

  \section{Bad News for matrix multiplication lower bounds}\label{badnews}

\begin{theorem}
\cite{MR3761737,MR2996880} It is essentially game   over for rank methods.
More precisely,   one cannot   prove bounds stronger than $\ur(T)\geq 6m$ for
$T\in \BC^m\ot \BC^m\ot \BC^m$ using rank methods.
\end{theorem}

 Remarkably this result was discovered 
essentially simultaneously by computer scientists and algebraic
geometers  with two completely
different proofs. 

I briefly explain the algebraic geometry proof as articulated by J. Buczy\'{n}ski (personal communication, see
\cite[\S 10.2]{MR3729273} for more detail):
a collection of $r$ distinct  points on a variety $X$ defines a zero dimensional smooth scheme
of length $r$, and a point on $\s_r(X)$ defines a 
zero dimensional smoothable scheme
of length $r$ (more precisely a point in the span of such). Determinantal methods detect zero dimensional schemes of length $r$.
Moreover, zero dimensional schemes of length $r$, even those supported at just one point, quickly fill the ambient
projective space. To rephrase (in what follows, $R$ is a zero dimensional scheme):
The secant variety is
$$\s_r(X):=\overline{ \bigcup  \{ \langle R\rangle\mid {\rm length}(R)=r,\ 
{\rm support}(R)\subset X,\ R:{\rm smoothable} \} }.
$$  
Define the   {\it cactus variety} \cite{MR3121848}: 
$$\kappa_r(X):=\overline{ \bigcup  \{ \langle R\rangle\mid {\rm length}(R)=r,\ 
{\rm support}(R)\subset X\} }
$$  

Determinantal equations are equations for the cactus variety and the cactus variety fills
the ambient space when $r$ is small ($6m$ for tensors).

\section{How to continue? Use more symmetry!}\label{bordersub}
So far, lower bounds were obtained by exploiting symmetry of the variety 
$\s_r(Seg(\BP A\times \BP B\times \BP C))$. But the point $\Mn$ also has symmetry. 
Write $A=U^*\ot V$, $B=V^*\ot W$, $C=W^*\ot U$. Then
$\Mn$ is $\Id_U\ot \Id_V\ot \Id_W$ re-ordered. Here $\Id_U: U\ra U$ is the identity map. Recall that the $GL(U)$-module
$U^*\ot U$ decomposes as $\fsl(U)\op \langle \Id_U\rangle$.  Thus

 $$G_{\Mn}\supset GL(U)\times GL(V)\times GL(W)=GL_n^{\times 3}\subset
 GL_{n^2}^{\times 3}.
 $$
 
 How to exploit this symmetry?
 
Given
$T\in A\ot B\ot C$, one has 
$\ur(T)\leq r$ if and only if  there exists a curve  $E_t\subset G(r,A\ot B\ot C)$
such that  

i) For $t\neq 0$, $E_t$ is spanned by $r$ rank one elements, and

ii) $T\in E_0$.
 
\medskip

Notice that if $E_t$ is such a curve, 
for all  $g\in G_T$, $gE_t$ also works.  This led to the following observation with M. Micha{\l}ek:

\begin{proposition} \label{normprop}\cite{MR3633766}
One can insist that  $E_0$ be fixed by a Borel subgroup of $G_T$. In particular,  for $\Mn$, 
one may insist that  $E_0$ is fixed by the action of triples of upper triangular $\nnn\times \nnn$ matrices
on $\La r((\BC^{\nnn^2} )^{\ot 3})$. 
\end{proposition}

This, combined with a border rank version of the classical
substitution method (see, e.g., \cite{MR3025382}),  led to what at the time I viewed as a Phyrric victory:

\begin{theorem}\cite{MR3842382}
    $\ur(\Mn)\geq 2 \nnn^2-\tlog_2\nnn-1$. 
    \end{theorem} 

I write Phyrric because it was clear this was the limit of the method.
Little did I realize that soon after, W. Buczy\'{n}ska and J. Buczy\'{n}ski would generalize
Proposition \ref{normprop} in a way that not only allowed further progress but also gave 
a potential path to overcoming the cactus barrier.

\section{Border Apolarity}\label{bapolar}
Buczy\'{n}ska-Buczy\'{n}ski had the following idea to use more information\cite{BBapolar}:   
Instead of considering limits of $r$ planes $\langle T_1(t)\hd T_r(t)\rangle$, where $T=\tlim_{t\ra 0}\sum T_j(t)$, consider     limits of {\it ideals} $I^t$, where    
$I^t$ is the  ideal of $[T_1(t)]\sqcup \cdots \sqcup [T_r(t)]$.    

This leads to the problem: How to take limits?  
A natural first idea is in the the {\it Hilbert scheme}. 

\subsection{The Hilbert Scheme of points}\label{hs}  One insists on saturated 
(with respect to the maximal ideal) ideals $I\subset Sym(V^*)$.    Then, in a sufficiently high degree $D$,
$I_D\subset S^DV^*$ determines $I$ in all degrees,   and \lq\lq sufficiently high\rq\rq\  can be made precise.   
Thus one is reduced to   to taking limits in one fixed Grassmannian. The Hilbert scheme
parametrizes saturated ideals with same {\it Hilbert polynomial}.

Let $I\subset Sym( V^*)$ be any ideal. Let $r_d=\tdim (S^dV^*/I_d)$,
so the Hilbert function is  $h_I(d):=r_d$.   
 Castelnuovo-Mumford regularity implies that  if
 one  fixes the Hilbert function, there exists an explicit $D=D(h_I)$ such that
$I_D$ determines $I_{D'}$ for all $D'>D$.     
Moreover $h_I(x)$ is a polynomial  when $x>D$, called the   Hilbert polynomial.

Bad news: The Hilbert scheme  doesn't work. Consider a toy case
of $3$ points in $\pp 2$:
$[1,0,0],[0,1,0],[1,-1,t]$   
$t\neq 0$, $(I^t)_1=0$    and 
$(I^t)_2=\langle x_3^2+t^2x_1x_2,\  x_3^2-tx_1x_3, \ x_1x_3+x_2x_3\rangle$   
But $(I^0)_1=\langle x_3\rangle$,    $(I^0)_2=\langle x_3^2, x_1x_3+x_2x_3\rangle$.   
The problem is that the ideal of the limiting scheme in a fixed degree  
is not the limit of spans and one  loses information 
important for border rank decomposition.

\subsection{The multigraded Hilbert scheme}\label{mghs}
The solution is to use the {\it Haiman-Sturmfels multi-graded Hilbert scheme} \cite{MR2073194}: 
Consider the product of Grasmannians
$$ 
 G(r_1,V^*)\times G(r_2,S^2V^*)\ctimes G(r_D,S^DV^*)
 $$    
 and map $I\mapsto ([I_1]\times [I_2]\ctimes [I_D])$.     For each $\BZ_{\geq 0}$-valued
 function $h$, get a (possibly empty) subscheme parametrizing all ideals 
 $I$ with {\it Hilbert function}  $h_I=h$.    
 This is
 rigged such that limit $I$ of ideals has same Hilbert function as ideals $I^t$.

Buczy\'{n}ska-Buczy\'{n}ski show that in border rank decompositions,   for $t>0$ 
one may assume the points are in general position    which leads to a constant
   Hilbert function
 as soon as is possible.

In the  tensor case, one has 
  more information because one has  curves of points on $Seg(\BP A\times \BP B\times \BP C)$.
One obtains  ideals in 
$Sym(A\op B\op C)^*=\bigoplus_{s,t,u} S^sA^*\ot S^tB^*\ot S^uC^*$,   
which is  $\BZ^{\op 3}$-graded.
This leads to a 
Hilbert function that  depends on three arguments:  $h_I(s,t,u):= \tdim ( S^sA^*\ot S^tB^*\ot S^uC^*/I_{s,t,u})$.   
 By the general position assumption,   $h_I(s,t,u)=\tmin\{ r, \tdim S^sA^*\ot S^tB^*\ot S^uC^*\}$.

 Instead of single curve $E_t\subset G(r,A\ot B\ot C)$ limiting to a Borel
 fixed point, for each
 $(i,j,k)$ one gets a curve in
 $G(r,S^iA^*\ot S^jB^*\ot S^kC^*)$, and  Buczy\'{n}ska-Buczy\'{n}ski 
  show that one may assume   that each  curve limits to  a Borel fixed point.   
 
 \subsection{Consequences}  The upshot is an  algorithm that either produces all normalized candidate $I^0$'s
 or proves border rank $>r$ as follows:

If $\ur(T)\leq r$, there exists  a  multi-graded ideal $I=I^0$  satisfying:
\begin{enumerate}
\item $I$ is contained in the annihilator of $T$. 
This condition says $I_{110}\subset T(C^*)^\perp$, $I_{101}\subset T(B^*)^\perp $, $I_{011}\subset T(A^*)^\perp$ and $I_{111}\subset T^\perp\subset A^*\ot B^*\ot C^*$.   
\item For all $(ijk)$ with $i+j+k>1$, $\tcodim I_{ijk}=r$.   
\item each $I_{ijk}$ is Borel-fixed.   
\item $I$ is an ideal, so the multiplication maps $I_{i-1,j,k}\ot A^*\op I_{i,j-1,k}\ot B^* \op I_{i,j,k-1}\ot C^*\ra  S^iA^*\ot S^jB^*\ot S^kC^*$
have image contained in $I_{ijk}$.  (These are rank conditions.)

\end{enumerate}

\subsection{Results}
Recall that Strassen proved   $\ur(\Mthree)\geq 14$,       Ottaviani 
and I  showed $\ur(\Mthree)\geq 15$, and   Micha{\l}ek and I     showed $\ur(\Mthree)\geq 16$.
Using border apolarity, with A. Conner and A. Harper, we showed: 
\begin{theorem}\cite{CHLapolar}
$\ur(\Mthree)\geq 17$.
\end{theorem}
The known upper bound is $20$ \cite{Smirnov13}.
Interestingly, we are able to construct candidate ideals for border rank $17$ and
we are currently attempting to determine if these ideals actually come from border rank
decompositions using deformation theory.

Recall that  so far only $\ur(\Mtwo)$ was known among nontrivial matrix multiplication tensors.  
Let $M_{\langle \aaa,\bbb,\ccc\rangle}\in \BC^{\aaa\bbb}\ot \BC^{\bbb\ccc}\ot \BC^{\ccc\aaa}$
denote the rectangular matrix multiplication tensor. 
Using border apolarity, we show
\begin{theorem}\cite{CHLapolar} $\ur(M_{\langle 223\rangle} )=10$ and  $\ur(M_{\langle 233\rangle} )=14$.
\end{theorem}
I also remark that the method gives a very short, computer free algebraic  proof that $\ur(\Mtwo)=7$.

All previous techniques for border rank lower bounds were useless when
one of the three vector spaces has dimension much larger than the other.
We also showed:

\begin{theorem} \cite{CHLapolar}
 For all $\nnn>25$,  $\ur(M_{\langle 2\nnn\nnn\rangle} )
\geq \nnn^2+ 1.32\nnn+1$ and for  
  all $\nnn>14$,  $\ur(M_{\langle 3\nnn\nnn\rangle} )
\geq \nnn^2+ 2\nnn$.  
\end{theorem}

Previously, only $\ur(M_{\langle 2\nnn\nnn\rangle} )
\geq \nnn^2+ 1$ and  $\ur(M_{\langle 3\nnn\nnn\rangle} )
\geq \nnn^2+2$ were known.  Notice that this also shows that border apolarity may be used
for sequences of tensors, not just fixed small tensors.

Currently we are working to strengthen the border apolarity algorithm, 
 to implement it more
efficiently, to take into account more geometric information, and to use deformation theory
to overcome the cactus barrier.

\section{Strassen's laser method and geometry}\label{lasersect}

In this last section I describe a program to utilize geometry to obtain upper bounds.

\subsection{Strassen's laser method and its barriers}
There was steady  progress upper bounding $\o$ from 1969 to 1988 culminating
in $\o<2.3755$ \cite{MR91i:68058}. All progress since 1984  has been obtained using methods from 
probability, statistical mechanics and information theory. 
 Given a tensor $T\in A\ot B\ot C$,
define its $k$-th Kronecker power  $T^{\boxtimes k}:=T^{\ot k}\in (A^{\ot k})\ot (B^{\ot k})\ot (C^{\ot k})$,
that is one takes its $k$-th tensor power and considers it as a $3$-way tensor instead of a $3k$-way tensor.
Since we are unable to upper bound the border  rank of the matrix multiplication
tensor directly, the idea of Strassen's laser method is to start with a tensor where we can
upper bound its border rank (e.g., a tensor of minimal border rank), take a large Kronecker power of
it, and then show  the Kronecker power degenerates to a large matrix multiplication tensor to get an upper bound
on the border rank of the large matrix multiplication tensor (namely the bound $\ur(T)^k$). 
Here one says $T$ degenerates to $T'$ if $T'\in \ol{GL(A)\times GL(B)\times GL(C)\cdot T}$. 
  Thus if the original tensor is of low cost (border rank)
and it produces a large matrix multiplication tensor (high value), it gives a good upper bound on $\o$.
 I emphasize that this is not done explicitly, just the existence of such a degeneration
is proved using methods from information theory pioneered by Shannon \cite{MR0026286}.
All the bounds since 1988 have been obtained by using a single tensor, the \lq\lq big\rq\rq\ Coppersmith-Winograd
tensor, which I'll denote $CW_q\in (\BC^{q+2})^{\ot 3}$ (there is one such for each $q$). Previous to that, the champion was
the \lq\lq little\rq\rq\ Coppersmith-Winograd tensor, which I'll denote $cw_q\in (\BC^{q+1})^{\ot 3}$.
I will not   describe the method here, see e.g., \cite{BCS,MR91i:68058,blaserbook} for expositions.

From 1988 until 2011 there was no progress whatsoever and starting 2011 there was incremental progress
leading up to the current record in \cite{doi:10.1137/1.9781611976465.32}.

In 2014, \cite{MR3388238} gave explanations for the halting progress, and showed there was
a limit to what one could prove with $CW_q$ (the limit is around $\o<2.3$).
Further explorations of  limits were made in \cite{MR3984617,DBLP:conf/innovations/AlmanW18,MR3984631}.

\begin{remark} An approach to upper bounds using the discrete Fourier
transform for finite groups was proposed in \cite{CU}. This approach yields similar bounds
to Strassen's laser method and faces similar barriers \cite{MR3631613,MR3891090}.
\end{remark}

A geometric explanation of the limits is given in \cite{MR3984631}:

Define 
the {\it asymptotic rank} of $T$:
$$\asrk(T):=\tlim_{N\ra\infty} (\ur(T^{\boxtimes  N}))^{\frac 1N},
$$
and the {\it asymptotic subrank} of $T$:
$$\assrk(T):=\tlim_{N\ra\infty} (\uq(T^{\boxtimes N}))^{\frac 1N}.
$$ 
For a given tensor $T$, a limit to its utility for the laser method is given by the ratio of these two quantities,
and the tensor could potentially be used to prove $\o$ is two only if the ratio is one, i.e., the
tensor is of minimal asymptotic rank and maximal asymptotic subrank. (In \cite{MR3984631}
they take the ratio of the logs, which they call {\it irreversibility}.) The barriers say nothing about
just how useful the tensor can be, only what one cannot do with it.

The {\it only} tensors we know the asymptotic rank of are those of minimal border rank.

The tensors $cw_q$ for $q<10$ could potentially be used to prove
$\o<2.3$, the main obstruction to doing so is that they are not of minimal border rank $q+1$
but instead $\ur(cw_q)=q+2$. Moreover, the case $cw_2$ could potentially be used to prove $\o$ is two.
What counts is the asymptotic rank, so  were
$\ur(cw_q^{\boxtimes 2})<\ur(cw_q)^2$, one could get a better upper bound on $\o$ than
the one found by Coppersmith-Winograd using the tensor.

Unfortunately for upper bounds, in \cite{CGLVkron} A. Conner, F. Gesmundo, E. Ventura and I  showed $\ur(cw_q^{\boxtimes 2})=\ur(cw_q)^2$
for $q>2$. At the time we were unable to determine the behaviour of $cw_2^{\boxtimes 2}$ as 
existing techniques did not yield any meaningful bound. With the advent of border apolarity, 
recently in \cite{CHLlaser}, A. Conner, H. Huang and I showed $\ur(cw_2^{\boxtimes 2})=\ur(cw_2)^2$.

I bring all this up in an article on representation theory and geometry because of recent work  in the search
for tensors useful for the laser method. Although the recent work has yet to improve upon the exponent
with geometric methods, it shows promise for the future.
I emphasize that in the computer science literature, the tensors used in the laser method
were found and exploited because of their combinatorial properties when expressed in a good
basis. M. Micha{\l}ek and I had the idea \cite{MR3682743,MR3633766}
to    analyze the  geometry of the tensors that have been successful
in proving upper bounds on $\o$ via the laser method, and then
to find other tensors with similar geometry in the hope they might be better for the laser method.
We found they had remarkable geometric properties. Among them, perhaps the most interesting
property 
is that  their   symmetry groups have large dimension.

\subsection{Tensors with symmetry}
There is a slight subtlety when discussing symmetry. The map
$GL(A)\times GL(B)\times GL(C)\ra GL(A\ot B\ot C)$ has a two dimensional kernel, namely
$\{(\lam \Id_A, \mu\Id_B, \nu\Id_C)\mid \lam\mu\nu=1\}$ and sometimes it is more convenient
to express the symmetry group (resp. algebra) in $GL(A)\times GL(B)\times GL(C)$ including
this kernel
(resp. $\fgl(A)\op \fgl(B)\op \fgl(C)$). When I do this I will decorate it with a tilde. 

First note that any minimal border rank tensor in $(\BC^m)^{\ot 3}$ has symmetry group of dimension
at least $2m-2$ as that is true for $\Mone^{\op m}:=a_1\ot b_1\ot c_1+\cdots + a_m\ot b_m\ot c_m$
and $\s_m(Seg(\pp {m-1}\times \pp{m-1}\times \pp{m-1}))=\ol {GL_m^{\times 3}\cdot [\Mone^{\op m}]}$.
A rank one tensor will have the largest symmetry group (of dimension $3m^2-3m+1$) but
tensors useful for the laser method have tended to be $1$-generic, so one expects a much 
smaller symmetry group.

As I illustrate below,  the tensor $CW_{m-2}$ has a symmetry group of dimension $\frac {m^2}2+\frac m2$,
which is quite large, so it is natural to look among $1$-generic tensors with large symmetry groups to
find ones useful for the laser method. This was the starting point of \cite{2019arXiv190909518C} with
A. Conner, F. Gesmundo, and E. Ventura.

Let $\bbta\in \langle a_2\hd a_{m-1}\rangle\ot \langle b_2\hd b_{m-1}\rangle\subset A\ot B$
be a nondegenerate bilinear form on $\BC^{m-2}\times \BC^{m-2}$.

\begin{theorem} \label{symthmc} \cite{2019arXiv190909518C}
Let $m \geq 7$ and let $\dim A = \dim B = \dim C = m$. Let $T \in A \otimes B \otimes C$ be a $1$-generic tensor. Then
\begin{equation}\label{eqn: bound for 1 generic tensor}
 \dim G_T < \frac{m^2}{2} + \frac{m}{2}-2
\end{equation}
except when $T$ is isomorphic to 
\be\label{Tform} S_{\bbta}:=     a_1 \otimes b_1 \otimes c_m + a_1 \otimes b_m \otimes c_1 + a_m \otimes b_1 \otimes c_1 +  
 \textsum_{\rho =2} ^{m-1} a_1 \otimes b_\rho \otimes c_\rho + \textsum_{\rho = 2}^{m-1} a_\rho \otimes b_1 \otimes c_\rho 
+ \bbta\otimes c_1,
\ene
 where $\bbta\in A\ot B$ is one of the four following rank $m-2$  bilinear forms

\begin{align}
&\label{skewcwB} \textsum_{\xi=2}^{p+1} a_\xi \otimes b_{\xi + p} - a_{\xi + p } 
\otimes b_\xi \ \ m=2p {\rm \ even}\ \   (T_{skewCW,m-2}) \\
&
\textsum_{\rho=2}^{m-1} a_\rho \otimes b_\rho \ \ {\rm all } \ m  \ \  (T_{CW,m-2}) \\
& a_{m-1} \otimes b_{m-1} + \textsum_{\xi = 2}^{p}  \bigl( a_\xi \otimes b_{\xi + p -1} - a_{\xi + p -1} \otimes b_\eta\bigr) 
 \ \ m=2p {\rm \ even} \ \   (T_{s+skewCW,m-2}) \\
&a_{m-1} \otimes b_{m-1} + \textsum_{\xi = 2}^{p}  \bigl( a_\xi \otimes b_{\xi + p -1} - a_{\xi + p -1} \otimes b_\eta\bigr) 
 \ \ m=2p+1 {\rm \ odd}\ \   (T_{s\oplus skewCW,m-2} )
\end{align}

All these except $T_{skewCW,m-2}$ have $\tdim G_T=\frac{m^2}{2}+\frac{m}{2}-2$, and
$\tdim G_{T_{skewCW,m-2}}=\frac{m^2}{2}+\frac {3m}{2}-4$.

In particular:  
when $m$ is even, there is a unique  up to isomorphism, $1$-generic tensor $T$ with maximal dimensional
symmetry group, namely $T_{skewCW,m-2}$, and there are 
exactly two, up to isomorphism, additional $1$-generic tensors $T$ such that $\dim G_T \geq \frac{m^2}{2}+\frac{m}{2}-2$, which are $T_{CW,m-2}$ and $T_{s+skewCW,m-2}$, where equality holds.

When $m$ is odd, there are exactly two  $1$-generic  tensors   $T$ up to isomorphism 
with maximal dimensional symmetry group $\frac{m^2}{2}+\frac{m}{2}-2$,  which are $T_{CW,m-2}$ and $T_{s\oplus skewCW,m-2}$.
\end{theorem}

\begin{remark}
In \cite{MR3754619} T. Seyannaeve decomposed $S^3\fgl_n$ and noticed that several of the highest
weight vectors that appeared were Coppersmith-Winograd tensors. This gave rise to the idea that
one might look among the highest weight vectors in $S^3\fgl_n$ to find ones useful for the laser
method. This was carried out in \cite{homs2021bounds}. This is a variant on having a large symmetry
group, as highest weight vectors are preserved by a parabolic subgroup.
\end{remark}

Call a tensor {\it skeletal} if it may be written in the form \eqref{Tform} for some nondegenerate
bilinear form $\bbta$.
 
\begin{proposition} \cite{2019arXiv190909518C}
\
 \begin{enumerate}
\item Any $1$-generic tensor in $(\BC^m)^{\ot 3}$ may be degenerated to
a skeletal tensor. 
\item The only skeletal tensor of minimal border rank is  
the Coppersmith-Winograd tensor, which is the case of 
 $\bbta\in S^2\BC^{m-2}$.
 
\item In particular any $1$-generic minimal border rank tensor $(\BC^m)^{\ot 3}$ degenerates to
$CW_{m-1}$.
\end{enumerate}
\end{proposition}

The result (3) originally appeared in \cite{hoyois2021hermitian}, although the proof (but not the statement)
was already in an early  preprint version of \cite{2019arXiv190909518C}.

In some sense (3) could be interpreted as saying that $CW_q$ is the worst minimal border rank $1$-generic tensor
for the laser method. The question is now, whether all others are equally bad, or just that the laser method
as currently practiced is not
refined enough.

I exhibit the symmetry Lie algebras of the above tensors:
Let $\tilde \fg_T$
denote the  Lie subalgebra of $ \fgl(A) \oplus \fgl(B) \oplus \fgl(C)  $ 
annihilating  $T$:
\[
\tilde \fg_T := \{ L \in  \fgl(A) \oplus \fgl(B) \oplus \fgl(C)  \mid L.T = 0\}.
\]
Here $L.T$ denotes the Lie algebra action.

If $L=(U,V,W)\in \fgl(A) \oplus \fgl(B) \oplus \fgl(C)$,
and we have bases $\{u^i_j\}, \{v^i_j\},\{w^i_j\}$ respectively for $U,V,W$, the condition $L.T = 0$ is equivalent to 
the following    system of linear  equations:
\begin{equation}\label{killeqn}
\sum_{i'}u^i_{i'} T^{i' jk}+\sum_{j'} v^j_{j'} T^{i j' k} + \sum_{k'}w^k_{k'} T^{ij k'}=0, \text{ for 
every 
} i,j,k. 
\end{equation}

For a skeletal $1$-generic tensor:
\begin{align}\label{LieTb}
&\tilde\fg_{S_{\bbta}}=
\\
&\nonumber
\left\{
\begin{pmatrix} u^1_1 & \ol u^\bt & u^1_m \\ 0 & X-\frac 12 u^1_1\Id  &   \bbta  \ol v - \ol z\\
0&0& -2u^1_1\end{pmatrix},
\begin{pmatrix} u^1_1 & \ol v^\bt & v^1_m \\ 0 & X-\frac 12 u^1_1\Id  &   \bbta\ol u -\ol z\\
0&0& -4u^1_1\end{pmatrix},
\begin{pmatrix} u^1_1 & \ol z^\bt   & -u^1_m-v^1_m \\ 0 & -X^{\bt}-\frac 12 u^1_1\Id  & -(\ol v+\ol u) \\
0&0& -4u^1_1\end{pmatrix} \right\}.
\end{align}
Here $\ol u,\ol v,\ol z\in \BC^{m-2}$ and $X\in \fh_\bbta$. 

 The small Coppersmith-Winograd tensor is 
$$cw_{m-2}=   +  
 \textsum_{\rho =2} ^{m-1} a_1 \otimes b_\rho \otimes c_\rho + \textsum_{\rho = 2}^{m-1} a_\rho \otimes b_1 \otimes c_\rho 
+ \bbta\otimes c_1 \in (\BC^{m-1})^{\ot 3}
$$
with $\bbta$ symmetric, and for other $\bbta$, write the tensors as $T_{\bbta-cw,m-2}$.
Then 
   $$
  \tilde \fg_{T_{\bbta-cw,q}}=
  \left\{ \left( \begin{pmatrix}-\mu-\nu & 0 \\ 0&\lambda\Id+ X \end{pmatrix},
  \begin{pmatrix}-\lambda-\nu & 0 \\ 0&\mu\Id+  X \end{pmatrix},
  \begin{pmatrix}-\lambda-\mu& 0 \\ 0&\nu\Id+  X \end{pmatrix}\right)
  \mid \lambda,\mu,\nu\in \BC
    X\in \fh_{\bbta} \right\}.
  $$
  In particular $\tdim \fg_{T_{cw,q}}=\binom{q}2+1$.

 All these tensors could potentially be used to prove $\o<2.3$ for $q<10$ and the case $q=2$ could
 again potentially be used to prove $\o$ is two. 
   Notice that if we take $\bbta$ to be completely skew, we get a tensor with a larger symmetry group.
 This tensor $T_{skewcw,2}$ unfortunately has larger initial cost than $cw_2$,
 namely $\ur(T_{skewcw,2})=5>4=\ur(cw_2)$,  (with the same value),
 however in \cite{CHLapolar}, using border apolarity,  we proved
 
 \begin{theorem} $\ur(T_{skewcw,2}^{\boxtimes 2})=17<< 25$\end{theorem}
 which provides hope for the laser method.
 
 \begin{remark} The Kronecker squares of $T_{skewcw,2}, cw_2$ are familiar tensors, respectively
 $\tdet_3$ and $\tperm_3$ considered as tensors, so these results are of interest well beyond
 the laser method.
 \end{remark}
 
 \section{Appendix: Strassen's alogorithm}\label{appen}

Here is Strassen's algorithm for multiplying 
$2\times 2$ matrices   using  $7$ scalar multiplications \cite{Strassen493}:
Set
\begin{align}
\label{starptwo}I&= (\aa 11 + \aa 22)(\bb 11 + \bb 22),\\ 
\nonumber II&=(\aa 21 + \aa 22)\bb 11, \\
\nonumber III&= \aa 11(\bb12-\bb 22)\\
\nonumber IV&=\aa 22(-\bb 11+\bb 21)\\
\nonumber V&=(\aa 11+\aa 12)\bb 22\\
\nonumber VI&= (-\aa 11+\aa 21)(\bb11 +\bb12),\\
\nonumber VII&=(\aa 12 -\aa 22)(\bb 21 + \bb 22),
\end{align}

\exerone{(1) Show that if $C=AB$, then
\begin{align*}
\cc 11&= I+IV-V+VII,\\
\cc 21&= II+IV,\\
\cc 12 &= III + V,\\
\cc 22 &= I+III-II+VI.
\end{align*}}

Now notice the the entries of $A,B$ themselves could be matrices, so this also gives,
by iterating, an algorithm for multiplying $2^k\times 2^k$ matrices.

\bibliographystyle{amsplain}

\bibliography{Lmatrix}

\def\cdprime{$''$} \def\cprime{$'$} \def\cprime{$'$} \def\cprime{$'$}
  \def\Dbar{\leavevmode\lower.6ex\hbox to 0pt{\hskip-.23ex \accent"16\hss}D}
  \def\cprime{$'$} \def\cprime{$'$} \def\cdprime{$''$} \def\cprime{$'$}
  \def\cprime{$'$} \def\Dbar{\leavevmode\lower.6ex\hbox to 0pt{\hskip-.23ex
  \accent"16\hss}D} \def\cprime{$'$} \def\cprime{$'$} \def\cprime{$'$}
  \def\cprime{$'$} \def\Dbar{\leavevmode\lower.6ex\hbox to 0pt{\hskip-.23ex
  \accent"16\hss}D} \def\cprime{$'$} \def\cprime{$'$}
\providecommand{\bysame}{\leavevmode\hbox to3em{\hrulefill}\thinspace}
\providecommand{\MR}{\relax\ifhmode\unskip\space\fi MR }
\providecommand{\MRhref}[2]{%
  \href{http://www.ams.org/mathscinet-getitem?mr=#1}{#2}
}
\providecommand{\href}[2]{#2}
\begin{thebibliography}{100}

\bibitem{MR2452824}
Hirotachi Abo, Giorgio Ottaviani, and Chris Peterson, \emph{Induction for
  secant varieties of {S}egre varieties}, Trans. Amer. Math. Soc. \textbf{361}
  (2009), no.~2, 767--792. \MR{MR2452824 (2010a:14088)}

\bibitem{AgrawalVinay}
M.~Agrawal and V.~Vinay, \emph{Arithmetic circuits: A chasm at depth four}, In
  Proc. 49th IEEE Symposium on Foundations of Computer Science (2008), 67–75.

\bibitem{AH}
J.~Alexander and A.~Hirschowitz, \emph{Polynomial interpolation in several
  variables}, J. Algebraic Geom. \textbf{4} (1995), no.~2, 201--222.
  \MR{96f:14065}

\bibitem{MR3025382}
Boris Alexeev, Michael~A. Forbes, and Jacob Tsimerman, \emph{Tensor rank: some
  lower and upper bounds}, 26th {A}nnual {IEEE} {C}onference on {C}omputational
  {C}omplexity, IEEE Computer Soc., Los Alamitos, CA, 2011, pp.~283--291.
  \MR{3025382}

\bibitem{MR3984617}
Josh Alman, \emph{Limits on the universal method for matrix multiplication},
  34th {C}omputational {C}omplexity {C}onference, LIPIcs. Leibniz Int. Proc.
  Inform., vol. 137, Schloss Dagstuhl. Leibniz-Zent. Inform., Wadern, 2019,
  pp.~Art. No. 12, 24. \MR{3984617}

\bibitem{doi:10.1137/1.9781611976465.32}
Josh Alman and Virginia~Vassilevska Williams, \emph{A refined laser method and
  faster matrix multiplication}, pp.~522--539.

\bibitem{DBLP:conf/innovations/AlmanW18}
\bysame, \emph{Further limitations of the known approaches for matrix
  multiplication}, 9th Innovations in Theoretical Computer Science Conference,
  {ITCS} 2018, January 11-14, 2018, Cambridge, MA, {USA}, 2018,
  pp.~25:1--25:15.

\bibitem{MR3388238}
Andris Ambainis, Yuval Filmus, and Fran{\c{c}}ois Le~Gall, \emph{Fast matrix
  multiplication: limitations of the {C}oppersmith-{W}inograd method (extended
  abstract)}, S{TOC}'15---{P}roceedings of the 2015 {ACM} {S}ymposium on
  {T}heory of {C}omputing, ACM, New York, 2015, pp.~585--593. \MR{3388238}

\bibitem{angelini2021description}
Elena Angelini and Luca Chiantini, \emph{On the description of identifiable
  quartics}, 2021.

\bibitem{MR3819100}
Elena Angelini, Luca Chiantini, and Nick Vannieuwenhoven, \emph{Identifiability
  beyond {K}ruskal's bound for symmetric tensors of degree 4}, Atti Accad. Naz.
  Lincei Rend. Lincei Mat. Appl. \textbf{29} (2018), no.~3, 465--485.
  \MR{3819100}

\bibitem{MR2500087}
Sanjeev Arora and Boaz Barak, \emph{Computational complexity}, Cambridge
  University Press, Cambridge, 2009, A modern approach. \MR{2500087
  (2010i:68001)}

\bibitem{MR1639346}
Sanjeev Arora, Carsten Lund, Rajeev Motwani, Madhu Sudan, and Mario Szegedy,
  \emph{Proof verification and the hardness of approximation problems}, J. ACM
  \textbf{45} (1998), no.~3, 501--555. \MR{1639346}

\bibitem{MR695915}
M.~D. Atkinson, \emph{Primitive spaces of matrices of bounded rank. {II}}, J.
  Austral. Math. Soc. Ser. A \textbf{34} (1983), no.~3, 306--315. \MR{MR695915
  (84h:15017)}

\bibitem{MR2996880}
Alessandra Bernardi and Kristian Ranestad, \emph{On the cactus rank of cubics
  forms}, J. Symbolic Comput. \textbf{50} (2013), 291--297. \MR{2996880}

\bibitem{MR605920}
D.~Bini, \emph{Relations between exact and approximate bilinear algorithms.
  {A}pplications}, Calcolo \textbf{17} (1980), no.~1, 87--97. \MR{605920
  (83f:68043b)}

\bibitem{MR592760}
Dario Bini, Grazia Lotti, and Francesco Romani, \emph{Approximate solutions for
  the bilinear form computational problem}, SIAM J. Comput. \textbf{9} (1980),
  no.~4, 692--697. \MR{MR592760 (82a:68065)}

\bibitem{blaserbook}
Markus Bl{\"a}ser, \emph{Fast matrix multiplication}, Graduate Surveys, no.~5,
  Theory of Computing Library, 2013.

\bibitem{MR3578455}
Markus Bl\"aser and Vladimir Lysikov, \emph{On degeneration of tensors and
  algebras}, 41st {I}nternational {S}ymposium on {M}athematical {F}oundations
  of {C}omputer {S}cience, LIPIcs. Leibniz Int. Proc. Inform., vol.~58, Schloss
  Dagstuhl. Leibniz-Zent. Inform., Wadern, 2016, pp.~Art. No. 19, 11.
  \MR{3578455}

\bibitem{MR3631613}
Jonah Blasiak, Thomas Church, Henry Cohn, Joshua~A. Grochow, Eric Naslund,
  William~F. Sawin, and Chris Umans, \emph{On cap sets and the group-theoretic
  approach to matrix multiplication}, Discrete Anal. (2017), Paper No. 3, 27.
  \MR{3631613}

\bibitem{MR3368091}
Grigoriy Blekherman and Zach Teitler, \emph{On maximum, typical and generic
  ranks}, Math. Ann. \textbf{362} (2015), no.~3-4, 1021--1031. \MR{3368091}

\bibitem{MR0660280}
Richard~P. Brent, \emph{The parallel evaluation of general arithmetic
  expressions}, J. Assoc. Comput. Mach. \textbf{21} (1974), 201--206.
  \MR{0660280 (58 \#31996)}

\bibitem{MR1243152}
Michel Brion, \emph{Stable properties of plethysm: on two conjectures of
  {F}oulkes}, Manuscripta Math. \textbf{80} (1993), no.~4, 347--371.
  \MR{MR1243152 (95c:20056)}

\bibitem{MR1601139}
\bysame, \emph{Sur certains modules gradu\'es associ\'es aux produits
  sym\'etriques}, Alg\`ebre non commutative, groupes quantiques et invariants
  ({R}eims, 1995), S\'emin. Congr., vol.~2, Soc. Math. France, Paris, 1997,
  pp.~157--183. \MR{1601139 (99e:20054)}

\bibitem{BBapolar}
Weronika Buczy\'{n}ska and Jaros{\l}aw Buczy\'{n}ski, \emph{Apolarity, border
  rank and multigraded {H}ilbert scheme}, arXiv:1910.01944, to appear in Duke
  Math J.

\bibitem{MR3121848}
\bysame, \emph{Secant varieties to high degree {V}eronese reembeddings,
  catalecticant matrices and smoothable {G}orenstein schemes}, J. Algebraic
  Geom. \textbf{23} (2014), no.~1, 63--90. \MR{3121848}

\bibitem{MR3769376}
Jaros{\l}aw Buczy\'{n}ski, Kangjin Han, Massimiliano Mella, and Zach Teitler,
  \emph{On the locus of points of high rank}, Eur. J. Math. \textbf{4} (2018),
  no.~1, 113--136. \MR{3769376}

\bibitem{BCS}
Peter B{\"u}rgisser, Michael Clausen, and M.~Amin Shokrollahi, \emph{Algebraic
  complexity theory}, Grundlehren der Mathematischen Wissenschaften
  [Fundamental Principles of Mathematical Sciences], vol. 315, Springer-Verlag,
  Berlin, 1997, With the collaboration of Thomas Lickteig. \MR{99c:68002}

\bibitem{Burgisser_2019}
Peter Burgisser, Cole Franks, Ankit Garg, Rafael Oliveira, Michael Walter, and
  Avi Wigderson, \emph{Towards a theory of non-commutative optimization:
  Geodesic 1st and 2nd order methods for moment maps and polytopes}, 2019 IEEE
  60th Annual Symposium on Foundations of Computer Science (FOCS) (2019).

\bibitem{MR3868002}
Peter B\"{u}rgisser, Christian Ikenmeyer, and Greta Panova, \emph{No occurrence
  obstructions in geometric complexity theory}, J. Amer. Math. Soc. \textbf{32}
  (2019), no.~1, 163--193. \MR{3868002}

\bibitem{MR2861717}
Peter B{\"u}rgisser, J.~M. Landsberg, Laurent Manivel, and Jerzy Weyman,
  \emph{An overview of mathematical issues arising in the geometric complexity
  theory approach to {${\rm VP}\neq{\rm VNP}$}}, SIAM J. Comput. \textbf{40}
  (2011), no.~4, 1179--1209. \MR{2861717}

\bibitem{MR2579394}
Dustin~A. Cartwright, Daniel Erman, Mauricio Velasco, and Bianca Viray,
  \emph{Hilbert schemes of 8 points}, Algebra Number Theory \textbf{3} (2009),
  no.~7, 763--795. \MR{2579394}

\bibitem{MR3023462}
Luca Chiantini and Giorgio Ottaviani, \emph{On generic identifiability of
  3-tensors of small rank}, SIAM J. Matrix Anal. Appl. \textbf{33} (2012),
  no.~3, 1018--1037. \MR{3023462}

\bibitem{MR3624400}
Luca Chiantini, Giorgio Ottaviani, and Nick Vannieuwenhoven, \emph{On generic
  identifiability of symmetric tensors of subgeneric rank}, Trans. Amer. Math.
  Soc. \textbf{369} (2017), no.~6, 4021--4042. \MR{3624400}

\bibitem{MR3984631}
Matthias Christandl, P\'{e}ter Vrana, and Jeroen Zuiddam, \emph{Barriers for
  fast matrix multiplication from irreversibility}, 34th {C}omputational
  {C}omplexity {C}onference, LIPIcs. Leibniz Int. Proc. Inform., vol. 137,
  Schloss Dagstuhl. Leibniz-Zent. Inform., Wadern, 2019, pp.~Art. No. 26, 17.
  \MR{3984631}

\bibitem{CU}
H~Cohn and C.~Umans, \emph{A group theoretic approach to fast matrix
  multiplication}, Proceedings of the 44th annual Symposium on Foundations of
  Computer Science (2003), no.~2, 438--449.

\bibitem{Como94:SP}
P.~Comon, \emph{Independent {C}omponent {A}nalysis, a new concept~?}, Signal
  Processing, Elsevier \textbf{36} (1994), no.~3, 287--314, Special issue on
  Higher-Order Statistics.

\bibitem{CGLVkron}
Austin Conner, Fulvio Gesmundo, J.M. Landsberg, and Emanuele Ventura,
  \emph{Kronecker powers of tensors and {S}trassen's laser method},
  arXiv:1909.04785.

\bibitem{2019arXiv190909518C}
Austin {Conner}, Fulvio {Gesmundo}, Joseph~M. {Landsberg}, and Emanuele
  {Ventura}, \emph{{Tensors with maximal symmetries}}, arXiv e-prints (2019),
  arXiv:1909.09518.

\bibitem{CHLapolar}
Austin Conner, Alicia Harper, and J.M. Landsberg, \emph{New lower bounds for
  matrix mulitplication and $\tdet_3$}, arXiv:1911.07981.

\bibitem{CHLlaser}
Austin Conner, Hang Huang, and J.M. Landsberg, \emph{Bad and good news for
  strassen's laser method: Border rank of $\tperm_3$ and strict
  submultiplicativity}, arXiv:2009.11391.

\bibitem{MR91i:68058}
Don Coppersmith and Shmuel Winograd, \emph{Matrix multiplication via arithmetic
  progressions}, J. Symbolic Comput. \textbf{9} (1990), no.~3, 251--280.
  \MR{91i:68058}

\bibitem{MR2996363}
Harm Derksen, \emph{Kruskal's uniqueness inequality is sharp}, Linear Algebra
  Appl. \textbf{438} (2013), no.~2, 708--712. \MR{2996363}

\bibitem{MR3761737}
Klim Efremenko, Ankit Garg, Rafael Oliveira, and Avi Wigderson, \emph{Barriers
  for rank methods in arithmetic complexity}, 9th {I}nnovations in
  {T}heoretical {C}omputer {S}cience, LIPIcs. Leibniz Int. Proc. Inform.,
  vol.~94, Schloss Dagstuhl. Leibniz-Zent. Inform., Wadern, 2018, pp.~Art. No.
  1, 19. \MR{3761737}

\bibitem{MR3126552}
Michael~A. Forbes and Amir Shpilka, \emph{Explicit {N}oether normalization for
  simultaneous conjugation via polynomial identity testing}, Approximation,
  randomization, and combinatorial optimization, Lecture Notes in Comput. Sci.,
  vol. 8096, Springer, Heidelberg, 2013, pp.~527--542. \MR{3126552}

\bibitem{MR2996364}
Shmuel Friedland, \emph{On tensors of border rank {$l$} in {$\Bbb{C}^{m\times
  n\times l}$}}, Linear Algebra Appl. \textbf{438} (2013), no.~2, 713--737.
  \MR{2996364}

\bibitem{MR1685640}
William Fulton, \emph{Eigenvalues of sums of {H}ermitian matrices (after {A}.
  {K}lyachko)}, no. 252, 1998, S\'{e}minaire Bourbaki. Vol. 1997/98, pp.~Exp.
  No. 845, 5, 255--269. \MR{1685640}

\bibitem{garfinkle}
D.~Garfinkle, \emph{A new construction of the {J}oseph ideal}, PhD thesis, MIT,
  1982.

\bibitem{MR0132079}
Murray Gerstenhaber, \emph{On dominance and varieties of commuting matrices},
  Ann. of Math. (2) \textbf{73} (1961), 324--348. \MR{0132079 (24 \#A1926)}

\bibitem{Gre11}
Bruno Grenet, \emph{{An Upper Bound for the Permanent versus Determinant
  Problem}}, Theory of Computing (2014), Accepted.

\bibitem{guo2021variety}
Zeyu Guo, \emph{Variety evasive subspace families}, 2021.

\bibitem{DBLP:journals/eccc/GuptaKKS13}
Ankit Gupta, Pritish Kamath, Neeraj Kayal, and Ramprasad Saptharishi,
  \emph{Arithmetic circuits: A chasm at depth three}, Electronic Colloquium on
  Computational Complexity (ECCC) \textbf{20} (2013), 26.

\bibitem{MR1504330}
J.~Hadamard, \emph{Sur les conditions de d\'ecomposition des formes}, Bull.
  Soc. Math. France \textbf{27} (1899), 34--47. \MR{1504330}

\bibitem{MR2073194}
Mark Haiman and Bernd Sturmfels, \emph{Multigraded {H}ilbert schemes}, J.
  Algebraic Geom. \textbf{13} (2004), no.~4, 725--769. \MR{2073194}

\bibitem{MR0384816}
Robin Hartshorne, \emph{Varieties of small codimension in projective space},
  Bull. Amer. Math. Soc. \textbf{80} (1974), 1017--1032. \MR{MR0384816 (52
  \#5688)}

\bibitem{MR3171099}
Jonathan~D. Hauenstein, Christian Ikenmeyer, and J.~M. Landsberg,
  \emph{Equations for lower bounds on border rank}, Exp. Math. \textbf{22}
  (2013), no.~4, 372--383. \MR{3171099}

\bibitem{homs2021bounds}
Roser Homs, Joachim Jelisiejew, Mateusz Michałek, and Tim Seynnaeve,
  \emph{Bounds on complexity of matrix multiplication away from cw tensors},
  2021.

\bibitem{MR983608}
Roger Howe, \emph{{$({\rm GL}_n,{\rm GL}_m)$}-duality and symmetric plethysm},
  Proc. Indian Acad. Sci. Math. Sci. \textbf{97} (1987), no.~1-3, 85--109
  (1988). \MR{983608}

\bibitem{hoyois2021hermitian}
Marc Hoyois, Joachim Jelisiejew, Denis Nardin, and Maria Yakerson,
  \emph{Hermitian k-theory via oriented gorenstein algebras}, 2021.

\bibitem{MR3695867}
Christian Ikenmeyer and Greta Panova, \emph{Rectangular {K}ronecker
  coefficients and plethysms in geometric complexity theory}, Adv. Math.
  \textbf{319} (2017), 40--66. \MR{3695867}

\bibitem{MR1957568}
Russell Impagliazzo, \emph{Hardness as randomness: a survey of universal
  derandomization}, Proceedings of the {I}nternational {C}ongress of
  {M}athematicians, {V}ol. {III} ({B}eijing, 2002), Higher Ed. Press, Beijing,
  2002, pp.~659--672. \MR{1957568}

\bibitem{MR3999690}
Joachim Jelisiejew, \emph{Elementary components of {H}ilbert schemes of
  points}, J. Lond. Math. Soc. (2) \textbf{100} (2019), no.~1, 249--272.
  \MR{3999690}

\bibitem{jelisiejew2021components}
Joachim Jelisiejew and Klemen Šivic, \emph{Components and singularities of
  quot schemes and varieties of commuting matrices}, 2021.

\bibitem{MR1663641}
Hajime Kaji, \emph{Secant varieties of adjoint varieties}, vol.~14, 1998,
  Algebra Meeting (Rio de Janeiro, 1996), pp.~75--87. \MR{1663641}

\bibitem{MR555701}
George Kempf and Linda Ness, \emph{The length of vectors in representation
  spaces}, Algebraic geometry ({P}roc. {S}ummer {M}eeting, {U}niv.
  {C}openhagen, {C}openhagen, 1978), Lecture Notes in Math., vol. 732,
  Springer, Berlin, 1979, pp.~233--243. \MR{555701}

\bibitem{koirand4}
Pascal Koiran, \emph{Arithmetic circuits: the chasm at depth four gets wider},
  preprint arXiv:1006.4700.

\bibitem{MR0444690}
Joseph~B. Kruskal, \emph{Three-way arrays: rank and uniqueness of trilinear
  decompositions, with application to arithmetic complexity and statistics},
  Linear Algebra and Appl. \textbf{18} (1977), no.~2, 95--138. \MR{MR0444690
  (56 \#3040)}

\bibitem{MR1923198}
Shrawan Kumar, \emph{Kac-{M}oody groups, their flag varieties and
  representation theory}, Progress in Mathematics, vol. 204, Birkh\"auser
  Boston Inc., Boston, MA, 2002. \MR{MR1923198 (2003k:22022)}

\bibitem{MR3314828}
\bysame, \emph{A study of the representations supported by the orbit closure of
  the determinant}, Compos. Math. \textbf{151} (2015), no.~2, 292--312.
  \MR{3314828}

\bibitem{MR2188132}
J.~M. Landsberg, \emph{The border rank of the multiplication of {$2\times2$}
  matrices is seven}, J. Amer. Math. Soc. \textbf{19} (2006), no.~2, 447--459.
  \MR{2188132 (2006j:68034)}

\bibitem{MR2865915}
\bysame, \emph{Tensors: geometry and applications}, Graduate Studies in
  Mathematics, vol. 128, American Mathematical Society, Providence, RI, 2012.
  \MR{2865915}

\bibitem{MR3729273}
\bysame, \emph{Geometry and complexity theory}, Cambridge Studies in Advanced
  Mathematics, vol. 169, Cambridge University Press, Cambridge, 2017.
  \MR{3729273}

\bibitem{MR1890196}
J.~M. Landsberg and Laurent Manivel, \emph{Construction and classification of
  complex simple {L}ie algebras via projective geometry}, Selecta Math. (N.S.)
  \textbf{8} (2002), no.~1, 137--159. \MR{1 890 196}

\bibitem{LMsec}
\bysame, \emph{On the ideals of secant varieties of {S}egre varieties}, Found.
  Comput. Math. \textbf{4} (2004), no.~4, 397--422. \MR{MR2097214
  (2005m:14101)}

\bibitem{MR3682743}
J.~M. Landsberg and Mateusz Micha{\l}ek, \emph{Abelian tensors}, J. Math. Pures
  Appl. (9) \textbf{108} (2017), no.~3, 333--371. \MR{3682743}

\bibitem{MR3633766}
\bysame, \emph{On the geometry of border rank decompositions for matrix
  multiplication and other tensors with symmetry}, SIAM J. Appl. Algebra Geom.
  \textbf{1} (2017), no.~1, 2--19. \MR{3633766}

\bibitem{MR3081636}
J.~M. Landsberg and Giorgio Ottaviani, \emph{Equations for secant varieties of
  {V}eronese and other varieties}, Ann. Mat. Pura Appl. (4) \textbf{192}
  (2013), no.~4, 569--606. \MR{3081636}

\bibitem{MR3724217}
J.~M. Landsberg and Nicolas Ressayre, \emph{Permanent v. determinant: an
  exponential lower bound assuming symmetry and a potential path towards
  {V}aliant's conjecture}, Differential Geom. Appl. \textbf{55} (2017),
  146--166. \MR{3724217}

\bibitem{MR3842382}
Joseph~M. Landsberg and Mateusz Micha{\l}ek, \emph{A {$2n^2-\log_2(n)-1$} lower
  bound for the border rank of matrix multiplication}, Int. Math. Res. Not.
  IMRN (2018), no.~15, 4722--4733. \MR{3842382}

\bibitem{MR3376667}
Joseph~M. Landsberg and Giorgio Ottaviani, \emph{New lower bounds for the
  border rank of matrix multiplication}, Theory Comput. \textbf{11} (2015),
  285--298. \MR{3376667}

\bibitem{LVdV}
R.~Lazarsfeld and A.~Van~de Ven, \emph{Topics in the geometry of projective
  space}, DMV Seminar, vol.~4, Birkh\"auser Verlag, Basel, 1984, Recent work of
  F. L. Zak, With an addendum by Zak. \MR{MR808175 (87e:14045)}

\bibitem{MR86c:68040}
Thomas Lickteig, \emph{A note on border rank}, Inform. Process. Lett.
  \textbf{18} (1984), no.~3, 173--178. \MR{86c:68040}

\bibitem{MR87f:15017}
\bysame, \emph{Typical tensorial rank}, Linear Algebra Appl. \textbf{69}
  (1985), 95--120. \MR{87f:15017}

\bibitem{MS8}
Ketan~D. Mulmuley, \emph{Geometric complexity theory: On canonical bases for
  the nonstandard quantum groups}, preprint.

\bibitem{MS6}
\bysame, \emph{Geometric complexity theory {VI}: the flip via saturated and
  positive integer programming in representation theory and algebraic
  geometry,}, Technical Report TR-2007-04, computer science department, The
  University of Chicago, May, 2007.

\bibitem{MS7}
\bysame, \emph{Geometric complexity theory {VII}: Nonstandard quantum group for
  the plethysm problem}, preprint.

\bibitem{MS5}
Ketan~D. Mulmuley and H.~Narayaran, \emph{Geometric complexity theory {V}: On
  deciding nonvanishing of a generalized {L}ittlewood-{R}ichardson
  coefficient}, Technical Report TR-2007-05, computer science department, The
  University of Chicago, May, 2007.

\bibitem{MS3}
Ketan~D. Mulmuley and Milind Sohoni, \emph{Geometric complexity theory {III}:
  on deciding positivity of {L}ittlewood-{R}ichardson coefficients}, preprint
  cs.ArXiv preprint cs.CC/0501076.

\bibitem{MS4}
\bysame, \emph{Geometric complexity theory {IV}: quantum group for the
  {K}ronecker problem}, preprint available at UC cs dept. homepage.

\bibitem{MS1}
\bysame, \emph{Geometric complexity theory. {I}. {A}n approach to the {P} vs.\
  {NP} and related problems}, SIAM J. Comput. \textbf{31} (2001), no.~2,
  496--526 (electronic). \MR{MR1861288 (2003a:68047)}

\bibitem{MS2}
\bysame, \emph{Geometric complexity theory. {II}. {T}owards explicit
  obstructions for embeddings among class varieties}, SIAM J. Comput.
  \textbf{38} (2008), no.~3, 1175--1206. \MR{MR2421083}

\bibitem{MR1344216}
David Mumford, \emph{Algebraic geometry. {I}}, Classics in Mathematics,
  Springer-Verlag, Berlin, 1995, Complex projective varieties, Reprint of the
  1976 edition. \MR{1344216 (96d:14001)}

\bibitem{nash2016open}
J.F. Nash and M.T. Rassias, \emph{Open problems in mathematics}, Springer
  International Publishing, 2016.

\bibitem{MR419491}
C.~Procesi, \emph{The invariant theory of {$n\times n$} matrices}, Advances in
  Math. \textbf{19} (1976), no.~3, 306--381. \MR{419491}

\bibitem{MR2719753}
Ran Raz, \emph{Elusive functions and lower bounds for arithmetic circuits},
  Theory Comput. \textbf{6} (2010), 135--177. \MR{2719753}

\bibitem{MR0506414}
Ju.~P. Razmyslov, \emph{Identities with trace in full matrix algebras over a
  field of characteristic zero}, Izv. Akad. Nauk SSSR Ser. Mat. \textbf{38}
  (1974), 723--756. \MR{0506414}

\bibitem{MR3891090}
Will Sawin, \emph{Bounds for matchings in nonabelian groups}, Electron. J.
  Combin. \textbf{25} (2018), no.~4, Paper No. 4.23, 21. \MR{3891090}

\bibitem{MR3754619}
Tim Seynnaeve, \emph{Plethysm and fast matrix multiplication}, C. R. Math.
  Acad. Sci. Paris \textbf{356} (2018), no.~1, 52--55. \MR{3754619}

\bibitem{MR0026286}
C.~E. Shannon, \emph{A mathematical theory of communication}, Bell System Tech.
  J. \textbf{27} (1948), 379--423, 623--656. \MR{MR0026286 (10,133e)}

\bibitem{Sipser}
Michael Sipser, \emph{The history and status of the p versus np question}, STOC
  '92 Proceedings of the twenty-fourth annual ACM symposium on Theory of
  computing (1992), 603--618.

\bibitem{Smirnov13}
A.V. Smirnov, \emph{The bilinear complexity and practical algorithms for matrix
  multiplication}, Computational Mathematics and Mathematical Physics
  \textbf{53} (2013), no.~12, 1781--1795 (English).

\bibitem{Strassen505}
V.~Strassen, \emph{Rank and optimal computation of generic tensors}, Linear
  Algebra Appl. \textbf{52/53} (1983), 645--685. \MR{85b:15039}

\bibitem{Strassen493}
Volker Strassen, \emph{Gaussian elimination is not optimal}, Numer. Math.
  \textbf{13} (1969), 354--356. \MR{40 \#2223}

\bibitem{MR3303254}
S{\'e}bastien Tavenas, \emph{Improved bounds for reduction to depth 4 and depth
  3}, Inform. and Comput. \textbf{240} (2015), 2--11. \MR{3303254}

\bibitem{MR763733}
B.~A. Trakhtenbrot, \emph{A survey of {R}ussian approaches to perebor
  (brute-force search) algorithms}, Ann. Hist. Comput. \textbf{6} (1984),
  no.~4, 384--400. \MR{763733}

\bibitem{vali:79-3}
Leslie~G. Valiant, \emph{Completeness classes in algebra}, Proc.~11th ACM STOC,
  1979, pp.~249--261.

\bibitem{Zak}
F.~L. Zak, \emph{Tangents and secants of algebraic varieties}, Translations of
  Mathematical Monographs, vol. 127, American Mathematical Society, Providence,
  RI, 1993, Translated from the Russian manuscript by the author.
  \MR{94i:14053}

\end{thebibliography}

\end{document}